\documentclass{amsart}

\newtheorem{theorem}{Theorem}[section]
\theoremstyle{plain}
\newtheorem{corollary}[theorem]{Corollary}
\newtheorem{lemma}[theorem]{Lemma}
\newtheorem{proposition}[theorem]{Proposition}
\theoremstyle{remark}

\newtheorem{example}[theorem]{Example}
\numberwithin{equation}{section}
\DeclareMathOperator{\tr}{tr}
\DeclareMathOperator{\otr}{0-tr}
\DeclareMathOperator{\ovol}{0-vol}
\DeclareMathOperator{\FP}{FP}
\DeclareMathOperator{\Area}{Area}
\DeclareMathOperator{\re}{Re}

\newcommand{\bbR}{\mathbb{R}}
\newcommand{\bbC}{\mathbb{C}}
\newcommand{\bbN}{\mathbb{N}}
\newcommand{\cinf}{C^\infty}
\newcommand{\del}{\partial}

\newcommand{\bX}{{\partial\bar X}}

\begin{document}
\title[Determinants and Isopolar Metrics]{Determinants of Laplacians and Isopolar Metrics on 
Surfaces of Infinite Area}
\author[Borthwick]{David Borthwick}
\address[Borthwick]{ Department of Mathematics and Computer Science, Emory University,
Atlanta, Georgia, 30322, U. S. A.}
\email{davidb@mathcs.emory.edu}
\author[Judge]{Chris Judge}
\address[Judge]{ Department of Mathematics, Indiana University, Bloomington, Indiana
47401, U. S. A.}
\email{cjudge@indiana.edu}
\author[Perry]{Peter A. Perry}
\address[Perry]{ Department of Mathematics, University of Kentucky, Lexington,
Kentucky, 40506--0027, U. S. A. }
\thanks{Borthwick supported in part by an NSF Postdoctoral Fellowship.}
\thanks{Judge supported in part by NSF\ grant DMS-9972425.}
\thanks{Perry supported in part by NSF\ grants DMS-9707051 and DMS-0100829.}
\date{September, 2001}

\begin{abstract}
We construct a determinant of the Laplacian for infinite-area surfaces
which are hyperbolic near infinity and without cusps.
In the case of a convex co-compact hyperbolic metric, the determinant can be related
to the Selberg zeta function and thus shown to be an entire function
of order two with zeros at the eigenvalues and resonances of the Laplacian.  
In the hyperbolic near infinity case the determinant is analyzed through the
zeta-regularized relative determinant for a conformal metric perturbation.  We establish
that this relative determinant is a ratio of entire functions of order two with divisor 
corresponding to 
eigenvalues and resonances of the perturbed and unperturbed metrics.  These results are 
applied to the problem of compactness in the smooth topology for the class of metrics with 
a given set of eigenvalues and resonances.
\end{abstract}
\maketitle
\tableofcontents

\section{Introduction}

Determinants of Laplacians provide an important, non-local spectral invariant
which plays a key role in the analysis of isospectral sets of compact
surfaces.  The goal of this paper is to develop the determinant of the
Laplacian for complete non-compact hyperbolic surfaces and their generalizations, 
as a tool in the spectral theory of these spaces. 
Such surfaces have infinite metric volume and at most finitely
many $L^{2}$-eigenvalues, so that `most' of the geometric information is
contained in the resonances (poles of the meromorphically continued resolvent).  
Thus it is natural to consider \textit{isopolar} classes of manifolds, those
with the same eigenvalues and resonances, and to investigate the restrictions imposed on 
their geometry.

For compact and finite-volume surfaces determinants are well-understood objects 
(see, for example \cite{DP:1986, 
Efrat:1988, OPS:1988a,OPS:1988b,Sarnak:1987,Voros:1988}).  In the infinite-volume
case there are several important points of contrast.  First and foremost,
it is not feasible to define the determinant through zeta-function regularization.  
Instead, we adopt a `Green's function method,' i.e. we define the determinant 
through the formal identity
\[ 
\left(\frac{d}{dz}\right)^{2}\log \det(A + z) =  - \tr (A+z)^{-2},
\]
valid for finite matrices $A$.  By regularizing the
trace of the square of the resolvent, one can produce a function $D(s)$, 
formally equal to $\det(\Delta + s(s-1))$.

The second issue is that the determinant is not a priori a `spectral' invariant; 
Indeed, a substantial part of our work will be to show that the eigenvalues and resonances 
do determine $D(s)$ up to finitely many parameters. 

A final major distinction to the compact case is that $\det(\Delta)$ so
defined is a proper function only on certain subsets of the moduli space of the underlying
topological surface (compare \cite{Khuri:1991} where a similar phenomenon is
analyzed for compact surfaces with boundary).  This lack of properness over the 
entire moduli limits the determinant's usefulness in determining whether or not 
isopolar classes are compact.

Throughout this paper we will consider only hyperbolic surfaces that are 
\textit{convex co-compact}, which means complete, finite topological type, 
infinite area, and without cusps. 
For such hyperbolic metrics, the results of Patterson-Perry \cite{PP:2000} allow us
to evaluate $D(s)$ in terms of Selberg's zeta function $Z_\tau(s)$ (our Theorem
\ref{thm.cc}), which is an entire function of order two.  
We thereby deduce that $D(s)$ has real zeros in
the half-plane $\re(s)>1/2$ corresponding to eigenvalues and zeroes corresponding
to resonances in the half-plane $\re(s)<1/2$ and possibly the special value 
$\re(s)=1/2$.   

Let $\bar X$ be a compact manifold with boundary, and $\rho$ be a defining coordinate
for the boundary $\bX$, i.e. $\rho\ge 0$, $\bX = \{\rho=0\}$, and $d\rho|_\bX \ne 0$.  
A \textit{conformally compact} metric is a metric on $X$ the form $g = \rho^{-2} \bar g$, 
where $\bar g$ is a smooth metric on $\bar X$.  The metric is called 
\textit{asymptotically hyperbolic} if $|dx|_{\bar g}|_{\bX} = 1$, which implies that 
Gaussian curvature of $g$ approaches $-1$ at $\bX$.  
(Note that this definition precludes cusps; the hyperbolic metrics which are conformally compact are
precisely the convex co-compact ones.)  Asymptotically hyperbolic is the level of generality
of the Mazzeo-Melrose \cite{MM:1987} parametrix construction and proof of meromorphic
continuation of the resolvent, although some extension to the conformally
compact case is possible \cite{B:2001}.
For an asymptotically hyperbolic metric one can still define the
determinant through the  resolvent as above, but there is no analogous zeta function theory for this
case.   Indeed, it is possible that the resolvent will have infinite-rank poles at
$-\frac12 \bbN_0$,  in which case $D(s)$ would fail to be entire.  
To avoid this, as well as
other  complications, we specialize to metrics that are \textit{hyperbolic near infinity}
(always assumed without cusps), meaning of constant curvature $-1$ outside some compact set
(same context as \cite{GZ:1995b}).

For metrics hyperbolic near infinity, in addition to the determinant defined as above we
can also construct a relative determinant based on conformal metric perturbation.
Recent results of Mazzeo-Taylor \cite{MT:2001}
imply that such a metric can be `uniformized' as $g = e^{2\varphi} \tau$,
where $\tau$ is a convex co-compact hyperbolic metric on $X$ and $\varphi 
\in \rho^2 \cinf(\bar X)$.   The decay of $\varphi$ at the boundary
allows the definition of a relative determinant of the Laplacians
$\Delta_g$ and $\Delta_\tau$, by zeta function 
regularization as in M\"uller \cite{Mueller:1998}.  The heat kernel definition of the
relative zeta function is
$$
\zeta(w,s)=\dfrac{1}{\Gamma(w)}\int_{0}^{\infty}t^{w}e^{-ts(s-1)}%
\tr \left[  e^{-t\hat\Delta_{g}}-e^{-t\Delta_\tau}\right]  \,\frac{dt}{t},
$$
where $\hat\Delta_g$ is the pull-back of $\Delta_g$ to $L^2(X,d\tau)$, and our
convention is that the Laplacian be a positive operator.  This converges for 
$\re(w)>1$ and $\re (s(s-1)) > -\mu$, where $\mu>0$ is the joint infimum of the spectra of 
the two Laplacians.
By analytic continuation to $w=0$ we may define
$$
D_{g,\tau}(s) = - \frac{d}{dw}\zeta_w(w,s)\Bigr|_{w=0},
$$
for $\re (s(s-1)) > -\mu$.

One would expect the relative determinant to have zeroes at the eigenvalues and
scattering poles of $g$, and poles at those of $\tau$.   In fact we can fully characterize the
relative determinant as a meromorphic function.   Let $P_{g}(s)$ and
$P_\tau(s)$ be Hadamard products formed from the eigenvalues and resonances, 
respectively, of $g$ and $\tau$ (see (\ref{eq.p0.def}) for a precise definition). 

\begin{theorem}
\label{thm.det.r}
Suppose $g$ is a metric hyperbolic near infinity, uniformized as $g = e^{2\varphi}\tau$.  
The relative determinant $D_{g,\tau}(s)$ extends to a meromorphic function of the form
\[
D_{g,\tau}(s)=e^{q(s)} P_g(s)/P_{\tau}(s).
\]
The polynomial $q(s)$ has degree at most two and is determined by the
eigenvalues and resonances of $\Delta_{g}$.
\end{theorem}

The proof of Theorem \ref{thm.det.r} consists of three steps. First, we
compute the divisor of $D_{g,\tau}(s)$ as follows.  We determine its zeros in
$\re(s)\geq1/2$ directly and obtain a functional equation of the form
\[
e^{h(s)}\frac{D_{g,\tau}(s)}{D_{g,\tau}(1-s)}=\det( S_{g,\tau}(s)),
\]
where $S_{g,\tau}(s)$ is the relative scattering operator, 
and $h(s)$ is a polynomial of degree at most two.  We
then to use results of Guillop\'e-Zworski \cite{GZ:1997} on determinants of scattering operators to
compute the divisor of the meromorphic function $\det(S_{g,\tau}(s))$.  
Secondly, we show that $D_{g,\tau}(s)$ is a
quotient of entire functions of order at most four by using estimates on the
relative zeta function together with constructive estimates on the resolvent
$R_g(s)$ proved in \cite{GZ:1997} (this is another point
where the restriction to hyperbolic near infinity is required).  Thirdly, we check that
$D_{g,\tau}(s)$ is entire of order two and prove the statement about
$q(s)$ by studying the asymptotics of $\log D_{g,\tau}(s)$ as $\re(s)
\rightarrow\infty$.

In the compact case there is a `Polyakov' formula expressing the relative determinant
for a conformal perturbation in terms of the conformal parameter, due to
Polyakov \cite{P:1981} and Alvarez \cite{A:1983}. 
Since the proof is based on the zeta regularization, the extension to our 
situation is quite straightforward.

\begin{proposition}
\label{prop.polyakov}
Suppose $g$ is an asymptotically hyperbolic metric with $K(g) + 1 = \mathcal{O}(\rho^2)$, 
uniformized as $g = e^{2\varphi}\tau$.  Then 
\[
\log D_{g,\tau}(1)=-\frac{1}{6\pi}\int_{X}\left(\frac12 |\nabla_{\tau}
\varphi|^{2}-\varphi\right)\,d\tau.
\]
\end{proposition}

With these tools in place, we turn to the isopolar problem.  Consider a 
topological surface $X$ of signature
$(h,M)$, i.e. $X$ is a surface with genus $h$ having $M$ discs removed. The
diffeomorphism classes are determined by the signature, and the Euler
characteristic is given by
\[
\chi(X)=2-2h-M.
\]
Let $\tau$ be a hyperbolic metric on $X$ which makes $X$ a complete Riemannian
manifold whose ideal boundary consists of $M$ circles. The surface $(X,\tau)$
takes the form
\[
\hat{X}\sqcup F_{1}\sqcup\cdots\sqcup F_{M}.
\]
Here $\hat{X}$, the convex core of $(X,\tau)$, is a convex compact manifold of
genus $h$ with geodesic boundary consisting of $M$ closed geodesics which we
denote by $\gamma_{1},\cdots,\gamma_{M}$. Letting $\ell_{i}$ be the geodesic
length of $\gamma_{i}$, the $F_{i}$ are hyperbolic funnels isometric to the
half-cylinder $(0,\infty)_{r}\times S_{\theta}^{1}$ with metric
\begin{equation}\label{fun.metric}
ds^{2}=dr^{2}+\ell_{i}^{2}\cosh^{2}r\,\,d\theta^{2}.
\end{equation}
The $F_{i}$ are glued to $\hat{X}$ along the bounding geodesics.

By studying the asymptotics of $D_{g,\tau}(s)$ as $\re(s)\to\infty$, we shall
prove the following:

\begin{proposition}\label{fix.euler}
Let $g$ be a metric hyperbolic near infinity on $X$.
The scattering poles and eigenvalues of $g$ determine the Euler characteristic of $X$. 
Thus the set of all surfaces $(X,g)$ with given eigenvalues and
scattering poles contains at most finitely many diffeomorphism types.
\end{proposition}

In fact, the eigenvalues and scattering poles of $g$ also determine the relative heat 
invariants for the pair $(g,\tau)$, which will be defined in \S\ref{sec.det.r}. 

In view of Proposition \ref{fix.euler}, we will fix a diffeomorphism type $(h,M)$ and a model
manifold $X$ with compactification $\bar X$.  
Suppose that $X$ carries a sequence of isopolar metrics $g_n$.
We wish to establish compactness of the isopolar set by showing that a 
subsequence of the $g_n$'s converges (modulo diffeomorphism) in a 
$\cinf$ topology.   Roughly the strategy is as follows:  uniformize $g_n = e^{2\varphi_n}
\tau_n$, and define relative determinants $D_{g_n,\tau_n}(s)$.  Note that
$D_{g_n,\tau_n}(s)$ is \emph{not} determined by the eigenvalues and resonances of $g_n$,
and therefore not independent of $n$  (in contrast to the 
Osgood-Phillips-Sarnak case \cite{OPS:1988b}).   
However, the relative heat invariants do turn out to be `spectral' invariants, 
and they give uniform estimates on $\varphi_n$ and its derivatives, expressed as integrals over $\tau_n$.  
These estimates may be combined with Theorem \ref{thm.det.r} and Proposition \ref{prop.polyakov},
to give uniform control over $Z_{\tau_n}(1)$.

We will demonstrate in Appendix \ref{sec.proper} that the evaluation of the Selberg 
zeta function at 1, $\tau\mapsto Z_{\tau}(1)$,  is not a proper function on the 
moduli space of $X$.  Let $\hat{X}_{\tau}$
denote the convex core of $(X,\tau)$ and let $\ell
(\del\hat{X}_{\tau})$ denote the sum of the lengths of the bounding geodesics.  
If we restrict to a subset of moduli space where $\ell(\del\hat{X}_{\tau})$ is bounded
above, then $Z_{\tau}(1)$ is a proper function.  With this restriction, then,
the uniform control of $Z_{\tau_n}(1)$ does imply convergence of a subsequence
of the $\tau_n$.  From this point on the strategy is identical to 
Osgood-Phillips-Sarnak, and the final result is:

\begin{theorem}
\label{thm.compact}
Let $g_n$ be a sequence of isopolar metrics of the form $g_n = e^{2\varphi_{n}} \tau_n$, where 
each $\varphi_n$ is supported in $\hat{X}_{\tau_n}$  and such that $\ell(\del\hat{X}_{\tau_n})$ is 
uniformly bounded in $n$.  Then there is a subsequence of the $g_{n}$ which converges, 
modulo diffeomorphisms of $X$, in the topology of $\rho^{-2}\cinf(\bar X; \mathcal{S}^2)$, to a
non-degenerate limiting metric in the same isopolar class.
\end{theorem}

Here $\mathcal{S}^2$ denotes the bundle of symmetric 2-tensors.  We note that
convergence in the weaker topology of $\cinf(X; \mathcal{S}^2)$ would not 
guarantee that the limit metric remains in the isopolar class.  

\medskip\noindent
\textbf{Remarks.}
\begin{enumerate}
\item
Examples of isopolar infinite-volume surfaces were given in Remark 2.3 of \cite{GZ:1997},
based on the transplantation method of B\'erard \cite{Be:1992}. 
Appendix D (contributed by Robert Brooks) presents some explicit cases that may be constructed 
by Sunada methods.
\item  
Our current methods require a priori information
concerning hyperbolic uniformization.  One might hope that $\ell(\del\hat{X}_{\tau_n})$
could be controlled uniformly by bounding the perimeter of the convex core of $g_n$, but we are
not aware of any such comparison a priori.  Our methods do not yield a comparison of this type
because all of the polar invariants are expressed in terms of the measures $d\tau_n$.  
Until convergence
of a subsequence of the $\tau_n$ is established, the relative heat invariants and the Polyakov
formula give no uniform information about the $\varphi_n$'s.  The uniform bound on
$\ell(\del\hat{X}_{\tau_n})$ must therefore be imposed explicitly to obtain the compactness of the
$\tau_n$ sequence.
\item  Assuming a convergent sequence of $\tau_n$'s, the invariants coming from the relative
determinant give uniform $H^m(X)$ bounds on $\varphi_n$ for all $m$.  This does not 
imply compactness of the sequence $\{\varphi_n\}$ in a $\cinf(\bar X)$ topology, which  
is the reason for  the restriction on the support of $\varphi_n$.   
\end{enumerate}

As a special case, one can consider isopolar classes of convex co-compact
hyperbolic metrics.  In Theorem \ref{thm.huber} we will prove an analog
of Huber's Theorem on the equivalence of the length spectrum and the
set of eigenvalues and resonances.  Then in Theorem \ref{ZFinite} we will show 
that the set of convex co-compact manifolds with the same length spectrum
is finite.  Together these imply:
\begin{theorem}
Let $R>0$.  Each set of isopolar convex co-compact hyperbolic surfaces 
with $\ell(\del\hat{X}_{\tau})< R$ is finite.
\end{theorem}

This theorem brings up an interesting question.  For a purely hyperbolic metric $\tau$,
can an upper bound for $\ell(\del\hat{X}_{\tau})$ be deduced from knowledge of
eigenvalues and resonances?  This is simpler than what would be required to generalize
the Theorem \ref{thm.compact}.
For a 1-holed torus ($h=1$, $M=1$), Buser-Semmler \cite{BS:1988}
show that there are no non-isometric surfaces with the same length spectrum, so
the answer is affirmative in the particular case.  No other cases appear to be known.

The results here complement the recent paper of Hassell-Zelditch
\cite{HZ:1999} where determinants of Laplacians on exterior planar domains in
Euclidean space are defined, and a compactness result for exterior domains
with the same scattering phase is proved. Roughly and informally, the
scattering phase is determined up to finitely many parameters by the
scattering poles, although this statement is difficult to make precise owing
to the lack of a sharp Poisson formula for resonances in this setting (see
\cite{Zworski:1998} for the best known results).

The plan of this paper is as follows.  In \S\ref{sec.spectral} we briefly 
review the spectral  and scattering theory for the
Laplacian on asymptotically hyperbolic surfaces, and then define
the determinant of the Laplacian.  This determinant is analyzed in the 
hyperbolic case in \S\ref{sec.det.constant}.  
In \S\ref{sec.det.r} we define and analyze the relative determinant for a conformal
metric perturbation, and prove Proposition \ref{prop.polyakov}.  
Theorem \ref{thm.det.r} is proved in \S\ref{sec.hadamard}. 
Finally, in \S\ref{sec.compact} we prove
the compactness theorems for isopolar metrics.  Appendix \ref{sec.proper}
contains the discussion of the properness of $Z_\tau(1)$ as a function on moduli space,
based on Theorem \ref{BersThm}, a generalization of Bers' theorem 
for pants decompositions of hyperbolic surfaces with geodesic boundary.
Appendices \ref{app.resolvent} and \ref{app.logdetsr} contain certain
technical facts needed elsewhere in the paper.  Finally, Appendix \ref{app.brooks},
contributed by Robert Brooks, discusses examples of isopolar surfaces arising from
the Sunada construction.

\textbf{Acknowledgment.}  Perry was partially supported by
a University Research Professorship from the University of Kentucky for the
academic year 1999-2000.  During a conference supported by the Research
Professorship, Lennie Friedlander set us straight about how to define determinants!
The work was completed in part at an MSRI workshop on Spectral Invariants
in May 2001, for which all three authors are grateful for support.

\section{Definition of the determinant}

\label{sec.spectral}

We begin by recalling some basic facts about the spectral theory of asymptotically 
hyperbolic surfaces.  For such a metric $g$, the Laplacian
$\Delta_g$ has at most finitely many eigenvalues in $[0,\frac{1}{4})$ 
(see \cite{LP:1982} for constant curvature and \cite{Mazzeo:1988, Mazzeo:1991b}
for variable curvature)
and absolutely continuous spectrum of infinite multiplicity in $[\frac{1}
{4},\infty)$ with no embedded eigenvalues (see e.g. \cite{LPAW} for constant
curvature and \cite{MM:1987} for the variable curvature). 
Thus the resolvent $(\Delta_g-z)^{-1}$ is a
meromorphic operator-valued function in the cut plane $\mathbb{C}
\backslash\lbrack\frac{1}{4},\infty)$.  Introducing the natural hyperbolic spectral
parameter $s$, we write 
$$
R_g(s)=(\Delta_g+s(s-1))^{-1}.
$$ 
Considered as a map from $L^{2}(X)$ to itself, $R_g(s)$
is then meromorphic in the half-plane $\re(s)>\frac{1}{2}$ with poles at real
numbers $\zeta>\frac{1}{2}$ for which $\zeta(1-\zeta)$ is an
eigenvalue of the Laplacian. 
Mazzeo and Melrose \cite{MM:1987} showed that, when
viewed as a map from $\cinf_0(X)$ to $\cinf(X)$, $R_g(s)$ admits a meromorphic 
extension to the complex plane.
Singularities of $R_g(s)$ with $\re(s)\leq\frac{1}{2}$ are called 
resonances (or scattering resonances).  We denote the full set of poles of $R_g(s)$  
(both eigenvalues and resonances) by $\mathcal{R}_g$.   A multiplicity can be assigned to each point
$\zeta\in\mathcal{R}_g$ as follows.  About each such $\zeta$,
the resolvent $R_g(s)$ has a Laurent expansion with finite polar part of the form
\[
\sum A_{j}(s-\zeta)^{-j}
\]
where the $A_{j}$ are finite-rank operators; we define the multiplicity,
$m_{\zeta}$, of a point $\zeta\in\mathcal{R}_g$ to be
\[
m_{\zeta}=\dim(\oplus_{j}\operatorname*{Ran}(A_{j})
).
\]
(For an $L^2$ eigenvalue this definition coincides with the usual notion of
multiplicity.)
For convenience, we will assume that poles are listed in $\mathcal{R}_g$ according to their
multiplicity.
In case $g$ has constant curvature $-1$, the existence of infinitely many
scattering resonances follows from Example 6, p. 856 and Remark 3, p. 851 of
\cite{SZ:1993}.

A `determinant of the Laplacian' should be an entire function
\[
D(s)=\det(\Delta_g+s(s-1))
\]
with zeros of multiplicity $m_{\zeta}$ at each point $\zeta\in\mathcal{R}_g$. To
motivate the definition we choose, consider the determinant
\[
D_{A}(s)=\det(A+s(s-1))
\]
where $A$ is a finite dimensional matrix.
One then has
\[
\left(\frac{1}{2s-1}\frac{d}{ds}\right)^{2}\log D_{A}(s)
= -\tr\left[(A+s(s-1))^{-2}\right]
\]
which suggests that $D(s)$ may be defined by replacing the right-hand side by
a suitable `trace' of $R_g(s)^{2}$.

Although $R_g(s)^{2}$ is not trace-class, its kernel is continuous.  Moreover,
$R_g(s)$ belongs to an algebra of pseudodifferential operators, the
$0$-pseudodifferential operators on $X$ (see \cite{MM:1987} or
\cite{Melrose:1995}), for which a natural renormalized trace is defined. 
To describe it, recall from the introduction that $X$ is assumed to be conformally
compact with respect to a boundary defining function $\rho$.
If $P$ is a 
$0$-pseudodifferential operator on $X$ with continuous kernel $\mathcal{K}_{P}$
(with respect to Riemannian measure), then
\begin{equation}\label{def.otr}
\otr (P)=\FP_{\varepsilon\downarrow0}\left(\int_{\rho
\geq\varepsilon}\mathcal{K}_{p}(x,x)\,dg(x)
\right)
\end{equation}
where $\FP_{\varepsilon\downarrow0}(\,\cdot\,)$ denotes the
Hadamard finite part.  (The structure of the $0$-calculus guarantees that the
argument has an asymptotic expansion in
$\varepsilon$ as $\varepsilon\to0$, and the Hadamard finite part is simply the constant term in
this expansion.)

In a similar way one can define the $0$-integral of any smooth function on $\bar X$,
and the $0$-volume of $X$ is just the $0$-integral of $1$.  Note that all of these definitions 
are dependent on the choice of $\rho$.  The $0$-integral of a smooth function depends on
the $1$-jet of $\rho$ restricted to $\bX$.

We will define a determinant $D_g(s)$ (up to two
free parameters plus dependence on the defining function) by the equation
\begin{equation}
\left(\frac{1}{2s-1}\frac{d}{ds}\right)^{2}\log D_g(s)=  
-\otr \left[R_g(s)^{2}\right]. \label{eq.ds}
\end{equation}
If $g$ is hyperbolic near infinity, then the 
poles of $R_g(s)$ are known to have finite-rank and one can hope that the function
$D_g(s)$ would be entire.   We will see that it is entire
provided that an appropriate defining function
is used to define the $0$-trace.  
If $g$ is only asymptotically hyperbolic then infinite-rank
poles at $- \frac12 \bbN_0$ cannot be ruled out.

\section{Properties of the determinant in the hyperbolic case}

\label{sec.det.constant}

In this section we will develop the theory of the determinant for the case
of a convex co-compact hyperbolic metric $\tau$ on $X$.
Recall that the Selberg zeta function $Z_\tau(s)$ is defined for $\re(s)>1$ as
a product over primitive closed geodesics $\gamma$ of $X$:
\begin{equation}
Z_\tau(s)=\prod_{\left\{  \gamma\right\}  }\prod_{k=0}^{\infty
}\left[  1- e^{-(s+k)\ell(\gamma)}\right]
\label{eq.zeta.euler}
\end{equation}
where $\ell(\gamma)$ is the length of $\gamma$. It is known
from \cite{PP:2000} (see Theorems 1.5 and 1.6) that $Z_\tau(s)$ is
an entire function of order two with zeros at the eigenvalues and scattering
resonances together the topological zeros of multiplicity 
$-(2k+1)\chi(X)$ at $s\in -\mathbb{N}_0$, where
$\chi(X)$ is the Euler characteristic of $X$.  Let $Z_{\infty}(s)$ be the
function (cf. Sarnak \cite{Sarnak:1987})
\begin{equation}
Z_{\infty}(s)=\left[\frac{(2\pi)^{s} \Gamma_{2}(s)^{2}}{\Gamma(s)}\right]^{-\chi(X)}. \label{eq.zinfty}
\end{equation}
Here $\Gamma_{2}(s)$ is Barnes' double gamma function \cite{Barnes:1900},
defined by the Hadamard product
\[
\frac{1}{\Gamma_{2}(s+1)}=(2\pi)^{s/2}
e^{-s/2-\frac{\gamma+1}{2}s^{2}}\prod_{k=1}^{\infty} \left(1+\frac{s}
{k}\right)^{k}e^{-s+s^{2}/k}
\]
with $\gamma$ Euler's constant. The function $Z_{\infty}(s)$
cancels the topological zeros of $Z_\tau(s)$, so that their product
is an entire function of order two. It is easily seen that
\begin{equation}
\frac{1}{\chi(X) (2s-1)}\frac{Z_{\infty}^{\prime}(
s)}{Z_{\infty}(s)}=-1-\frac{1}{s}-\gamma
+\sum_{k=1}^{\infty}\frac{s}{k(s+k)}. \label{eq.zinfty.d1}
\end{equation}

The Hadamard product for the eigenvalues and resonances is
\begin{equation}
P_\tau(s) = \prod_{\zeta\in\mathcal{R}_{\tau}} \left(1-\frac
{s}{\zeta}\right)^{m_{\zeta}}e^{m_{\zeta}(-\frac{s}{\zeta}+\frac{s^{2}}{2\zeta^{2}})}
\label{eq.ptau.def}.
\end{equation} 
By the characterization of $Z_\tau(s)$ in \cite{PP:2000}, we have
\begin{equation}\label{zeta.had}
Z_\tau(s) Z_\infty(s) = e^{q(s)} P_\tau(s),
\end{equation}
for some polynomial $q(s)$ of degree two.  In particular, the length spectrum
determines the eigenvalues and resonances.
Asymptotic analysis of (\ref{zeta.had}) leads to the following 
infinite-volume analog of Huber's Theorem:
\begin{theorem}\label{thm.huber}
The eigenvalues and resonances of $\Delta_\tau$ (with multiplicities)
determine both $\chi(X)$ and the function
$Z_\tau(s)$.   Thus for a convex co-compact hyperbolic metric the eigenvalues and resonances  
determine the length spectrum and vice versa.
\end{theorem}
\begin{proof}
The expansion of $\log Z_{\infty}(s)$ was used prominently by Sarnak
\cite{Sarnak:1987} and is based on classical results of Barnes
\cite{Barnes:1900} for the double gamma function. For $\re(s)\rightarrow
\infty$ we have
\begin{equation}
\begin{split}
\log Z_{\infty}(s)  &  \sim\chi(X)[\tfrac{1}{2}\log2\pi+\tfrac{1}{4}
-2\zeta^{\prime}(-1)-(\tfrac{1}{2}s(s-1)-\tfrac{1}{6})\log
s(s-1)\\
&  \qquad+\tfrac{3}{2}s(s-1)]+\sum_{l=1}^{\infty}c_{l}[s(s-1)]^{-l}.
\end{split}
\label{eq.zinf.exp}
\end{equation}
In particular, this expansion has a term of the form
\begin{equation}\label{log.term}
-\chi(X)(\tfrac12 s(s-1) - \tfrac16)\log s(s-1).
\end{equation}
It is clear from (\ref{eq.zeta.euler}) that $\log Z_\tau(s)
=\mathcal{O}(e^{-s\ell_{0}})$ as $\re(s)\rightarrow\infty$, where $\ell_{0}$ is
the length of the shortest closed geodesic on $X$.  Thus the left-hand side of (\ref{zeta.had})
has an asymptotic expansion with a term of the form (\ref{log.term}).  This could
not possibly be canceled by $q(s)$ on the right-hand side,  so $\chi(X)$ is determined
by $P_\tau(s)$ and hence by the set of eigenvalues and resonances.  

Once $\chi(X)$ is known, $q(s)$ is the only unknown in the asymptotic expansion of $\log P_\tau(s)$, 
so it too must be determined by the eigenvalues and resonances.  Note from above that $Z_\infty(s)$
depends only on $\chi(X)$.  Thus by (\ref{zeta.had}), $Z_\tau(s)$ is fixed by the
eigenvalues and resonances.

The proof that the length spectrum may be extracted from $Z_\tau(s)$ is completely analogous
to the compact case.  For example, one may define $\ell_0$ as the unique number $\omega$
such that for real $s$
$$
-\infty < \lim_{s\to\infty} e^{\omega s}\log Z_\tau(s) < 0.
$$
Then terms with $\ell_0$ are removed from the product, the same approach determines $\ell_1$,
etc.
\end{proof}

In order to connect the determinant to the zeta function we must 
be careful about the definition of the $0$-trace.  On a hyperbolic surface one
may specify a natural class of defining function by requiring that the $0$-volume
of the funnels $F_j$ equal zero.  For instance, in the model metric (\ref{fun.metric})
we may take $\rho = e^{-r}$.  If the $0$-volume of the funnels is zero, then
the $0$-volume of $X$ equals the volume of $\hat X_\tau$, which is $-2\pi \chi(X)$ 
by Gauss-Bonnet.  We will use $\otr_\tau$ to denote the $0$-trace with respect to
a defining function in this class.  This is the same convention used for the $0$-trace 
in \cite{GZ:1997, PP:2000}.

\begin{theorem}
\label{thm.cc}
Suppose that $\tau$ is a convex co-compact hyperbolic metric on $X$,
and let $D_\tau(s)$ be any function defined by (\ref{eq.ds}) using $\otr_\tau$.  
Then $D_\tau(s)$ is an entire function of order two with zeros $\zeta\in
\mathcal{R}_\tau$ of multiplicity $m_{\zeta}$, given by the formula
\begin{equation}
D_\tau(s)=e^{Fs(s-1)+G} Z_\tau(s)Z_{\infty}(s)\label{eq.ds.form}
\end{equation}
where $F$ and $G$ are the free parameters in the definition.  Therefore $D_\tau(s)$ is
determined (up to the constants $F$ and $G$) by the set $\mathcal{R}_\tau$ of eigenvalues and 
resonances of $\Delta_\tau$ counted with multiplicities.
\end{theorem}

\noindent\textbf{Remarks.}
\begin{enumerate}
\item Sarnak \cite{Sarnak:1987} showed that the (zeta-regularized) determinant of the
Laplacian on a compact hyperbolic surface $S$ is given by
\[
\det(\Delta_{S}+(s-1))=Z_{S}(s)  \left(e^{E-s(s-1)}(2\pi)^{s}\frac{\Gamma
_{2}(s)^{2}}{\Gamma(s)}\right)^{-\chi(S)}
\]
where $Z_{S}(s)$ is the zeta function for the compact surface,
and
\[
E=-\frac{1}{4}-\frac{1}{2}\log(2\pi)+2\zeta^{\prime}(-1),
\]
with $\zeta(s)$ the Riemann zeta function. In particular,
\begin{equation}
\det(\Delta_{S})=c(\chi(S))Z_{S}(1) \label{eq.det.surface}
\end{equation}
for a constant depending only on the Euler characteristic. Thus our determinant has the
same form as Sarnak's; from this point of view a natural choice for the
constants $F$ and $G$ in (\ref{eq.ds.form}) would be $F=\chi(X)$ and
$G=-\chi(X)E$. 
\item  Taking $s=1$ gives a determinant of
the Laplacian in terms of special values of the zeta function (in contrast to
\cite{Sarnak:1987} there is no derivative because neither the zeta function
nor the determinant has a zero at $s=1$): one has
\[
\det\Delta_{g}=e^{G}(2\pi)^{-\chi(X)}Z_\tau(1).
\]
\end{enumerate}

\begin{proof}
To prove Theorem \ref{thm.cc}, we study the function $L(s) = \frac{d}{ds} \log D_\tau(s)$,
defined up to one parameter by the relation
\[
\frac{1}{2s-1}\frac{d}{ds}\left(\frac{1}{2s-1}L(s)\right)
= - \otr_\tau R_\tau(s)^{2}.
\]
We recall from \cite{PP:2000} the relation
\begin{equation}\label{dlogz}
\frac{1}{2s-1}\frac{d}{ds}\log Z_\tau(s)=\otr_\tau
(R_\tau(s)-R_{\mathbb{H}}(s)),
\end{equation}
where $R_{\mathbb{H}}(s)$ is the resolvent of the Laplacian on
the two-dimensional hyperbolic space $\mathbb{H}$, and the integral kernel on
the right is obtained by lifting the integral kernel of $R_\tau(s)$ to
$\mathbb{H}\times\mathbb{H}$, subtracting that of $R_{\mathbb{H}}(
s)$, and restricting to the diagonal; the resulting restriction is
smooth and projects to a function on ${\bar X}$ (see \cite{PP:2000}, Section 6). 
The $0$-trace in (\ref{dlogz}) refers, by a slight abuse of notation,
to the $0$-integral of this smooth function.
A further differentiation gives
\begin{equation}\label{ddZ}
\left(\frac{1}{2s-1}\frac{d}{ds}\right)^{2}\log Z_\tau(s)
=\frac{1}{2s-1}\frac{d}{ds}\left(\frac{1}{2s-1}L(s)\right) 
+\otr_\tau(R_{\mathbb{H}}(s)^{2})
\end{equation}
where the second right-hand term is interpreted as follows.  The integral
kernel of $R_{\mathbb{H}}(s)^{2}$ is continuous and its
restriction to the diagonal of $\mathbb{H}{\times}\mathbb{H}$ equals the
constant
\[
-\frac{1}{2\pi}\frac{1}{2s-1}\frac{d}{ds}\psi(s),
\]
where $\psi(s)$ is the logarithmic derivative of $\Gamma(s)$. Recall for later use that
\begin{equation}
\psi(s)=-\frac{1}{s}-\gamma+\sum_{k=1}^{\infty}\frac
{s}{k(s+k)}. \label{eq.psi}
\end{equation}
As noted above, $\ovol_\tau(X)=-2\pi\chi(X)$, so we have
$$
\otr_\tau(R_{\mathbb{H}}(s)^{2}) = \frac{\chi(X)}{2s-1}\frac{d}{ds}\psi(s).
$$
From (\ref{ddZ}) we thus see that
\begin{equation}
L(s)=\frac{d}{ds}\log Z_\tau(s)+\chi(X)(
2s-1)\psi(s)+ C(2s-1)
\label{eq.L}
\end{equation}
for some constant $C$.

We now compare (\ref{eq.zinfty.d1}) and (\ref{eq.psi}) to conclude that
\[
\psi(s)=\frac{1}{\chi(X)}\frac{1}{2s-1}\frac{d}{ds}\log
Z_{\infty}(s)+1
\]
so that (\ref{eq.L}) becomes
\[
L(s)=\frac{d}{ds}\log(Z_\tau(s)Z_{\infty}(s))+(\chi(X)+C)(2s-1).
\]
Setting $F=C+\chi(X)$, we obtain the formula (\ref{eq.ds.form}).
\end{proof}

\section{Relative determinant}

\label{sec.det.r}

\subsection{Conformal deformation theory}
\label{sec.conformal}

On a compact surface any metric is conformal related to unique hyperbolic metric.
This is the basis for the isospectral result of Osgood-Phillips-Sarnak \cite{OPS:1988b}.
In the non-compact case,
Mazzeo-Taylor recently proved a uniformization
result well-suited to our situation.

\begin{theorem} \cite{MT:2001}
Let $h$ be an arbitrary smooth metric on $\bar X$.  Then there exists a function 
$u\in \rho C^\infty(\bar X)$ such that $e^{2u} \rho^{-2} h$ is a complete hyperbolic
metric on $X$.
\end{theorem}

(A theorem of Hulin-Troyanov \cite{HT:1992} gives a similar result in this situation
with $u$ bounded and smooth on $X$, but with no boundary regularity on $\bar X$.)
If $g$ denotes the conformally compact metric $\rho^{-2} h$, then $u$ is a solution 
of the equation
\begin{equation}\label{cdef.eq}
\Delta_g u = - K(g) - e^{2u}.
\end{equation}
By examining this closely we can gain an extra order of vanishing of $u$ with a 
curvature restriction.
\begin{corollary}\label{cor.rho2}
Let $g$ be a conformally compact metric on $X$.
If $K(g) + 1 = \mathcal{O}(\rho^2)$, then there exists a function 
$u \in \rho^2 C^\infty(\bar X)$
such that $e^{2u} g$ is hyperbolic on $X$.
\end{corollary}
\begin{proof}
The theorem gives $u\in \rho C^\infty(\bar X)$ satisfying (\ref{cdef.eq}), which we rewrite as
$$
(\Delta_g+2)u = -(1+K(g)) + 1 +2u - e^{2u}.
$$
The first term on the right-hand side is $\mathcal{O}(\rho^2)$ by assumption, and the remainder is
also because of the vanishing of $u$ to first order.  Thus $(\Delta_g+2)u = \mathcal{O}(\rho^2)$.
Now if we write $u = \rho f$, then because $g$ is an asymptotically hyperbolic metric on $X$,
a straightforward computation shows $(\Delta_g+2) (\rho f) = 2 \rho f + \mathcal{O}(\rho^2)$.
Thus $f$ vanishes at $\rho=0$ and $u$ is in fact in $\rho^2 C^\infty(\bar X)$.
\end{proof}

Note that the curvature condition is stronger than the asymptotically hyperbolic
assumption, which 
implies only that $K(g) + 1 = \mathcal{O}(\rho)$.

\subsection{Zeta-regularization of the relative determinant}

Let $g$ be an asymptotically hyperbolic metric on $X$ such that $K(g)+1 = \mathcal{O}(\rho^2)$.  
According to Corollary \ref{cor.rho2} this can be written as  
$g=e^{2\varphi} \tau$ where $\varphi\in \rho^2\cinf(\bar X)$ and $\tau$ is convex 
co-compact hyperbolic.   To define a relative determinant from $\Delta_{g}$
to $\Delta_{\tau}$ we need first to compare operators acting on the same space.
Note that $\Delta_g = e^{-2\varphi} \Delta_\tau$ and $dg = e^{2\varphi} d\tau$. 
Thus
\[
\mathcal{U}f=e^{\varphi}f
\]
is a unitary transformation from $L^{2}(X,dg)$ to $L^{2}(X,d\tau)$, and 
the pull-back of $\Delta_g$ under this map is
\[
\hat\Delta_g=\mathcal{U}\Delta_{g}\mathcal{U}^{-1}=e^{-\varphi}\Delta_{\tau}e^{-\varphi}.
\]

We will define the relative determinant through a relative zeta function, as in \cite{Mueller:1998}. 
For this to be well-defined we need $e^{-t\hat\Delta_g}-e^{-t\Delta_\tau}$ to be trace-class for $t>0$.
For $\varphi$ compactly supported this can be seen immediately from Duhamel's formula:
\[
e^{-t\hat\Delta_g}-e^{-t\Delta_{\tau}}=\int_{0}^{t}e^{-s\Delta_{\tau}}\left(  
\Delta_{\tau}-\hat\Delta_g\right)  e^{-\left(  t-s\right)  \Delta_{g}}ds.
\]
For the general case, we will prove a slightly stronger condition.  Let $R_\tau(s)$
be the resolvent for $\Delta_\tau$, and define 
$$
\hat R_g(s) = (\hat\Delta_g + s(s-1))^{-1} = e^{\varphi} R_g(s) e^{-\varphi}.
$$
Note that in the definition (\ref{eq.ds}) of $D_g(s)$, replacing
$R_g(s)$ by $\hat R_g(s)$ has no effect.
\begin{lemma}\label{lemma.restrace} 
Let $g$ be an asymptotically hyperbolic metric with $K(g)+1 = \mathcal{O}(\rho^2)$.  
For $\re(s)>1$ the operator $\hat R_g(s)^{2}-R_\tau(s)^{2}$ is trace-class.
\end{lemma}

\begin{proof}
A detailed description of the structure of the resolvent $R_g(s)$ was obtained in 
\cite{MM:1987}.  From this picture one can easily deduce that the operator
$\rho^\alpha \hat R_g(s) \rho^\beta$ is Hilbert-Schmidt provided that 
$\re(s)>1/2 - \min\{\alpha,\beta\}$ and $\alpha+\beta >1/2$.  The same
fact holds for $R_\tau(s)$.

We write
\begin{equation}\label{eq.res2.diff}\begin{split}
\hat R_g(s)^{2}-R_\tau(s)^{2}  &= \hat R_g(s) (\hat R_g(s) - R_\tau(s)) 
+ (\hat R_g(s) - R_\tau(s)) R_\tau(s) \\
&= \hat R_g(s)^2 (\Delta_\tau - \hat\Delta_g) R_\tau(s)
+ \hat R_g(s) (\Delta_\tau - \hat\Delta_g) R_\tau(s)^2\\
\end{split}\end{equation}
We will show that the first of these terms is a trace-class operator since the
analysis of the second term is similar.  Expanding
\[
\Delta_\tau - \hat\Delta_g=(1-e^{-\varphi})\Delta_\tau e^{-\varphi} + \Delta_\tau
(1- e^{-\varphi}),
\]
we can re-express the right-hand first term in (\ref{eq.res2.diff}) as
$$
\hat R_g(s)^{2}(e^{-\varphi}-1)\Delta_\tau e^{-\varphi} R_\tau(s)
 +\hat R_g(s)^{2}\Delta_\tau(e^{-\varphi}-1) R_\tau(s).
$$ 
Once again, the two terms are very similar and we will consider only the first.
Since $\Delta_\tau e^{-\varphi} R_\tau(s)$ is a bounded operator for $\re(s)>1/2$,
we may ignore this part.  The function $(e^{-\varphi}-1) = \rho^2 h$ for $h\in 
\cinf(\bar X)$, so the analysis reduces to considering $\hat R_g(s)^{2} \rho^2$.
Writing this as $\hat R_g(s)\rho^\beta \cdot \rho^{-\beta} \hat R_g(s) \rho^2$,
we see from the characterization of Hilbert-Schmidt operators above that this term 
is trace class for $\re(s) >1$.
\end{proof}

By the Birman-Krein spectral shift theory (which we will use in more detail below), 
we deduce the following:
\begin{corollary}
For $g$ an asymptotically hyperbolic metric with $K(g)+1 = \mathcal{O}(\rho^2)$,  
the relative heat-kernel $e^{-t\hat\Delta_g}-e^{-t\Delta_\tau}$ is trace-class for $t>0$.
\end{corollary}

As $t\rightarrow0$, standard heat kernel asymptotics can be used to derive an expansion
\begin{equation}
\tr \left[  e^{-t\hat\Delta_g}-e^{-t\Delta_\tau}\right]  \sim\frac{1}
{t}\sum_{j\geq0}a_{j}t^{j}. \label{eq.heatexp}
\end{equation}
Let $\mathcal{H}_{g}(t,x,y)$ be the heat kernel for the metric
$g$, i.e. the Schwarz kernel of $e^{-t\Delta_{g}}$ with respect to the
Riemannian measure $dg$. If $\hat\mathcal{H}_{g}(t,x,y)$
denotes the Schwarz kernel of $e^{-t\hat\Delta_g}$ with respect to $d\tau$, then the
relationship between these kernels is
\[
\hat\mathcal{H}_{g}(t,x,y)=e^{\varphi(x)}\mathcal{H}_{g}(t,x,y)e^{\varphi(y)}.
\]
Thus, the local heat expansion for $\hat\Delta_g$ is easily related to the standard
results for $\Delta_{g}$ (see \cite{MS:1967}). The zeroth relative heat
invariant measures the change in area:
\begin{equation}
a_{0}=\frac{1}{4\pi}\int_{X}(e^{2\varphi}-1)\,d\tau.
\label{eq.a0}
\end{equation}
The formula for the first invariant is
\begin{equation}
a_{1}=\frac{1}{12\pi}\int_{X}(e^{2\varphi}K_{g}-K_\tau )\,d\tau,
\label{eq.a1}
\end{equation}
where the $K$'s are the respective Gaussian curvatures.
Since $K_\tau =-1$ and $K_{g}=e^{-2\varphi}(\Delta_\tau\varphi
-1)$, we have
\[
a_{1}=\frac{1}{12\pi}\int_{X}\Delta_\tau\varphi\,d\tau=0.
\]
(The first heat invariant in the compact case is the Euler characteristic;  here the first relative 
heat invariant is zero because the Euler characteristic for $\tau$ to $g$ are the same.) 
The second heat invariant involves the square of the curvature:
\begin{equation}
a_{2}=\frac{1}{60\pi}\int_{X}(e^{2\varphi}K_{g}^{2}-1)\,d\tau,
\label{eq.a2}
\end{equation}
and the higher invariants are integrals of polynomials in $K_{g}$ and
$\Delta_{g}$:
\begin{equation}
a_{j}=c_{j}\int_{X}e^{2\varphi}K_{g}\Delta_{g}^{j-2}K_{g}\,d\tau+\text{(terms
with fewer derivatives)} \label{eq.aj}
\end{equation}
where $c_{j}\neq0$.

Let $\lambda_{g}=\inf\operatorname*{spec}(\hat\Delta_g)$ and $\lambda_{\tau}
=\inf\operatorname*{spec}(\Delta_\tau)$. Both of these numbers are strictly
positive since $X$ does not support constant eigenfunctions. If $\mu
=\min(\lambda_{g},\lambda_{\tau})$, then we can deduce that
\begin{equation}
\tr \left[  e^{-t\hat\Delta_g}-e^{-t\Delta_\tau}\right]  =\mathcal{O}(e^{-\mu t})
\label{eq.htrdecay}
\end{equation}
as $t\rightarrow\infty$.

The relative zeta function is initially defined for $\re(w)>1$ and
$\re(s(s-1))>-\mu$ by
\begin{equation}
\zeta(w,s)=\dfrac{1}{\Gamma(w)}\int_{0}^{\infty}t^{w}e^{-ts(s-1)}
\tr \left[  e^{-t\hat\Delta_g}-e^{-t\Delta_\tau}\right]  \,\frac{dt}{t}
\label{eq.zeta}
\end{equation}
We wish to define the relative determinant by
\begin{equation}
\log D_{g,\tau}(s)=-\zeta_{w}(0,s), \label{eq.rdet.def}
\end{equation}
where the subscript denotes differentiation with respect to the $w$ variable.
To justify this definition, note that the heat expansion (\ref{eq.heatexp})
implies that
\[
\zeta(w,s)=a_{0}\frac{[s(s-1)]^{w-1}}{w-1}+(\text{analytic for }\re(w)>-1),
\]
so (\ref{eq.rdet.def}) is well-defined by meromorphic continuation in $w$, 
at least for $\re(s(s-1))>-\mu$.

The next result will enable us to connect the relative and absolute
determinants.

\begin{lemma}
\label{lemma.dr.trace}The identity
\[
\left(\frac{1}{2s-1}\frac{d}{ds}\right)^{2}\log D_{g,\tau}
(s)=-\tr \left[ \hat R_g(s)^{2}-R_\tau(s)^{2}\right]
\]
holds for $\re(s)>1$.
\end{lemma}

\begin{proof}
For $\re\left[  s(s-1)\right] > -\mu$, the function $\zeta(w,s)$
is continuously differentiable to all orders in $s$ and $w$ near $w=0$ so that
we may calculate
\begin{align}
\left(\frac{1}{2s-1}\frac{d}{ds}\right)^{2}\log D_{g,\tau}(s)  &  =\left.
\frac{d}{dw} \Bigl[  -w(w+1)\zeta(w+2,s)\Bigr]\right|_{w=0} \label{eq.d2}\\
&  =-\zeta(2,s).\nonumber
\end{align}
On the other hand, $\hat R_g(s)^{2}-R_\tau(s)^{2}$ is trace-class for $\re(s)>1$ 
by Lemma \ref{lemma.restrace}.  For real $s$ with
$s>1$ and real $w$ with $w\geq2$ we can therefore identify
\begin{equation}\label{zeta.realw}
\zeta(w,s)=\tr \left[ \hat R_g(s)^{w}-R_\tau(s)^{w}\right]
\end{equation}
(the power is well-defined since the resolvents are positive operators if
$s>1$). The result now follows for all real $s>1$ and hence all $s$ with
$\re(s)>1$ by analytic continuation.
\end{proof}

In order to get the correct relation to the zeta function, in \S\ref{sec.det.constant} we defined
$D_\tau(s)$ using $\otr_\tau$, the $0$-trace for a class of defining function canonically 
associated to $\tau$.  Now, in order to get the proper connection to the relative determinant,
we must also use $\otr_\tau$ to define $D_g(s)$. 
Then, since
\[
\tr [ \hat  R_g(s)^{2}-R_\tau(s)^{2}]  =\otr_\tau [  \hat R_g(s)^{2}] -\otr_\tau [R_\tau(s)^{2}] 
\]
for $\re(s)>1$, it follows that
\[
\left(\frac{1}{2s-1}\frac{d}{ds}\right)^{2}\log D(s)= \left(\frac
{1}{2s-1}\frac{d}{ds}\right)^{2}\log(D_{g,\tau}(s)D_{\tau}(s)).
\]
Thus:
\begin{lemma}
\label{lemma.abs.and.rel}
Let $D_\tau(s)$ and $D_g(s)$ be defined by (\ref{eq.ds}) using $\otr_\tau$.
Then there are constants $E$ and $F$ so that
\[
D_g(s)=e^{E+Fs(s-1)}D_{g,\tau}(s)D_{\tau}(s)\text{.}
\]
In particular, if $g$ is hyperbolic near infinity then
$D_{g,\tau}(s)$ admits a meromorphic continuation to the whole plane.
\end{lemma}

To conclude this subsection, we note that the Birman-Krein theory of the
spectral shift (see e.g. \cite{Yafaev:1991}, chapter 8 for an exposition)
applies to give a measurable, locally integrable function $\xi$ on
$[0,\infty)$, the spectral shift function, with the property that
\begin{equation}
\int\phi^{\prime}(\lambda)d\xi(\lambda)
=\tr (\phi(\hat\Delta_g)-\phi(\Delta_{\tau}))\label{eq.krein}
\end{equation}
for any smooth function $\phi$ vanishing rapidly at infinity. In particular we
have
\[
\tr \left[  e^{-t\hat\Delta_g}-e^{-t\Delta_\tau}\right]  =-t
\int_{0}^{\infty}e^{-t\lambda}\xi(\lambda)\,d\lambda.
\]
It follows from (\ref{eq.htrdecay}) that the spectral shift function is
supported in $[\mu,\infty)$ where $\mu>0$. From Lemma \ref{lemma.restrace} we
obtain the estimate
\begin{equation}
\int_{\mu}^{\infty}\lambda^{-3}\left|  \xi(\lambda)\right|
\,d\lambda<\infty. \label{eq.xi.bd}
\end{equation}
This in turn allows us to derive the following formula:
\begin{equation}\label{xi.ss}
\xi(w,s)=-w\int_{\mu}^{\infty}(\lambda+s(s-1))^{-w}\,\xi(\lambda)\,d\lambda,
\end{equation}
which is valid for $\re(w)\geq2$ and $s(1-s)
\in\mathbb{C}\backslash\lbrack0,\infty)$.

\subsection{Polyakov formula}

\label{sec.polyakov}

In this subsection we will establish the Polyakov formula stated in the introduction.

\begin{proof}[Proof of Proposition \ref{prop.polyakov}]
Consider the variation $\varphi\rightarrow\varphi+\delta\varphi$ for some function
$\delta\varphi\in \rho^2 C^{\infty}(\bar X)$. 
Since $\hat\Delta_g=e^{-\varphi}\Delta_{\tau}e^{-\varphi}$ we have
\[
\delta \hat\Delta_g=-\delta\varphi\>\hat\Delta_g-\hat\Delta_g\delta\varphi.
\]
Because of the decay  of $\delta\varphi$ at the boundary, we can differentiate inside
the trace to obtain
\[
\begin{split}
\delta\tr\left[  e^{-t\hat\Delta_g}-e^{-t\Delta_\tau}\right]   &  =t\tr\left[
\delta\varphi\>\hat\Delta_g e^{-t\hat\Delta_g}+\hat\Delta_g\delta\varphi\>e^{-t\hat\Delta_g}\right] \\
&  =2t\tr\left[  \delta\varphi\>\hat\Delta_g e^{-t\hat\Delta_g}\right]  .
\end{split}
\]
Then the variation of the relative zeta function (at $s=1$) is
\[
\begin{split}
\delta\zeta(w,1)  &  =\frac{2}{\Gamma(w)}\int_{0}^{\infty}t^{w}\tr\left[
\delta\varphi\>\hat\Delta_g e^{-t\hat\Delta_g}\right]  \>dt\\
&  =-\frac{2}{\Gamma(w)}\int_{0}^{\infty}t^{w}\frac{d}{dt}\tr\left[
\delta\varphi\>e^{-t\hat\Delta_g}\right]  \>dt\\
&  =\frac{2w}{\Gamma(w)}\int_{0}^{\infty}t^{w-1}\tr\left[  \delta
\varphi\>e^{-t\hat\Delta_g}\right]  \>dt.
\end{split}
\]

Now from the expansion of the heat kernel at $t=0$ we have
\[
\tr\left[  \delta\varphi\>e^{-t\hat\Delta_g}\right]  =b_{0}t^{-1}+b_{1}+\mathcal{O}(t),
\]
where
\[
b_{1}=\frac{1}{12\pi}\int_{X}\delta\varphi\>(\Delta_\tau\varphi-1)\>d\tau.
\]
The heat expansion implies that
\[
\delta\zeta(w,1)=\frac{2w}{\Gamma(w)}\left[  \frac{b_{0}}{w-1}+\frac{b_{1}}
{w}+(\text{analytic for }\re(w)>-1)\right]  ,
\]
from which we directly obtain
\[
\delta\log D_{g,\tau}(1)=-2b_{1}=-\frac{1}{6\pi}\int_{X}\delta\varphi\>
(\Delta_{\tau}\varphi-1)\>d\tau.
\]
Integrating this expression gives the Polyakov formula.
\end{proof}

\section{Hadamard factorization of the relative determinant}

\label{sec.hadamard}

For this section and the remainder of the paper we will restrict our attention 
to the case of a metric $g$ which is hyperbolic near infinity.
so that $D_{g,\tau}(s)$ is a meromorphic function by Lemma \ref{lemma.abs.and.rel}.

Let us introduce the Hadamard product
\begin{equation}
P_{g}(s)=\prod_{\zeta\in\mathcal{R}_{g}} \left(1-\frac
{s}{\zeta}\right)^{m_{\zeta}}e^{m_{\zeta}(-\frac{s}{\zeta}+\frac{s^{2}}{2\zeta^{2}})}
\label{eq.p0.def},
\end{equation}
whose zeroes are given by $\mathcal{R}_g$.  Convergence of the Hadamard product is 
guaranteed by the estimate
\begin{equation}
\#\left\{  \zeta\in\mathcal{R}_{*}\text{\thinspace:\thinspace}\left|
\zeta\right|  \leq r\right\}  \leq C(1+r^{2}) \label{eq.r0.est}
\end{equation}
proven in \cite{GZ:1995a}. 

The goal of this section is the proof of Theorem \ref{thm.det.r}.
This will be divided into three steps: first, we show that the
divisors of $D_{g,\tau}(s)$ and the quotient $P_g(s)/P_{\tau}(s)$ coincide, so that the
two differ by an holomorphic function without zeros. Next, we show that the
logarithmic derivative of
\[
W(s)=D_{g,\tau}(s)\frac{P_{\tau}(s)}{P_g(s)}
\]
has at most polynomial growth.
Finally, we obtain an exact growth rate by examining the asymptotic behavior of $\log
D_{g,\tau}(s)$ as $s\rightarrow\infty$ through real values.

\subsection{Divisor of the relative determinant}

\label{sec.hadamard.divisor}

In this subsection, we prove:

\begin{proposition}
\label{prop.divisor}
For $g$ hyperbolic near infinity, the relative determinant  $D_{g,\tau}(s)$
has a meromorphic extension to $\bbC$, with divisor equal to that of $P_g(s)/P_{\tau}(s)$.
\end{proposition}

The connection between $D_{g,\tau}(s)$ and the resonances will be made through
scattering theory. Let $S_{g}(s)$ and $S_{\tau}(s)$ be respective 
scattering operators (on $\del\bar X$) for $\hat\Delta_g$ and $\Delta_\tau$ as defined in \cite{GZ:1997},
Definition 2.12.  By \cite{JS:1998}, Theorem 7.1 we have that $S_{g}(s)-S_{\tau}(s)$ 
is a pseudodifferential operator of order $2\re(s)-3$. Thus the
relative scattering operator $S_{g,\tau}(s)=S_{g}(s)S_{\tau}(s)^{-1}$ satisfies
$S_{g,\tau}(s)=I+Q(s)$, where $Q(s)$ is of order $-2$ and hence trace class.

Our first result follows from the analysis in Proposition 2.5 of \cite{GZ:1997}

\begin{lemma}
\label{sreldet} The determinant of the relative scattering operator satisfies
\[
\det S_{g,\tau}(s)=e^{f_1(s)}\frac{P_g(1-s)}{P_g(s)}\frac{P_{\tau}(s)}{P_{\tau}(1-s)},
\]
for some polynomial $f_1$ of degree at most four.
\end{lemma}

Since both $P_g(1-s)$ and $P_g(s)$ appear in Lemma \ref{sreldet}, in order to make
use of it we must first explicitly count poles of $D_{g,\tau}(s)$ in the half-plane
$\re(s) \geq1/2$.

\begin{lemma}
\label{divdr} The divisor of $D_{g,\tau}(s)$ in the closed half-plane $\re(s)
\ge 1/2$ coincides with that of $P_g(s)/P_{\tau}(s)$.
\end{lemma}

\begin{proof}
To determine the zeros and poles of $D_{g,\tau}(s)$ in the closed half-plane
$\re(s)\geq1/2$, we want to consider the trace of $\hat R_g(s)-R_\tau(s)$. This does
not even have a 0-trace, so instead we fix $s_{0}$ not a pole of either
operator and examine
\begin{equation}
H(s)=(2s-1)\otr_\tau \Bigl[ \hat R_g(s) - R_\tau(s) - \hat R_g(s_0) + R_\tau(s_0) \Bigr]. \label{eq.h.def}
\end{equation}
By Lemma {\ref{lemma.dr.trace}} we have
\[
\frac{1}{2s-1}\frac{d}{ds}\left(\frac{1}{2s-1}H(s)\right)=\left( \frac
{1}{2s-1}\frac{d}{ds}\right)^{2}\log D_{g,\tau}(s).
\]
Integrating gives
\begin{equation}
\frac{d}{ds}\log D_{g,\tau}(s)=H(s)+c(2s-1) \label{eq.rdet.dlog}
\end{equation}
for some constant $c$.  For $* = g$ or $\tau$, let $R_*(s)$ denote $\hat R_g(s)$ or $R_\tau(s)$,
respectively.  Let $m_\zeta^*$ be the multiplicity of 
the pole of $R_*(s)$ at $\zeta$, with $m_\zeta^*=0$ if there is no pole.
For $\re(\zeta)>\frac{1}{2}$, the only poles occur at real $\zeta$ and correspond to 
eigenvalues.  At such a value the pole of $R_*(s)$ is simple, and the
residue of $(2s-1)R_*(s)$ is a rank-$m_{\zeta}^*$
projection operator. 
From this it follows that for $\re(\zeta)>1/2$ $H(s)$ has a simple pole with residue exactly
$m_{\zeta}^g-m_{\zeta}^{\tau}$.  Then by (\ref{eq.rdet.dlog}) 
$D_{g,\tau}(s)$ has divisor $\{(\zeta,m_{\zeta}^g-m_{\zeta}^{\tau})\}$ for $\re(\zeta)>1/2$, 
which is the divisor of $P_g(s)/P_{\tau}(s)$.

The only singularity of $\hat R_g(s)$ or $R_\tau(s)$ that can occur on the line
$\re(s)=1/2$ is at $s=1/2$ (see \cite{GZ:1997}, Lemma 4.1). It follows from
that Lemma and the fact that there are no $L^{2}(X)$-eigenfunctions of
$\Delta_{g}$ or $\Delta_\tau$ with eigenvalue $1/4$ that
\[
R_*(s)=\frac{B_*(s)}{2s-1}+C_*(s),
\]
where $B_{*}(s)$ and $C_*(s)$ are analytic at $s=1/2$ and $B_*(s)$
is a finite-rank operator.  Although $(2s-1)R_*(s)$ is thus analytic near $s=1/2$, 
the definition of the zero-trace together
with the spatial decay of the kernels may introduce a pole there, as may be
seen by considering the small-$\varepsilon$ expansion of the model integral
$\int_{\varepsilon}^{1}y^{2s-2}dy$.   From \cite{PP:2000}, 
Lemma 4.9 we obtain that
\[
B_*(s;z,z')=\sum_{k=1}^{N_{\ast}}u_{k}^{\ast}(s;z) v_{k}^{\ast}(s;z'),
\]
with $u_{k}^{\ast}(s), v_{k}^{\ast}(s) \in \rho^{s} \cinf(\bar X)$.  These functions
have the property that the restrictions
$\bar{u}_{k}^{\ast}  = (\rho^{-1/2}u_{k}^{\ast}(1/2))|_{\bX}$
and $\bar{v}_{k}^{\ast}  = (\rho^{-1/2}v_{k}^{\ast}(1/2))|_{\bX}$ are both
non-zero.  Analyzing as in the proof of Theorem 6.2 of \cite{PP:2000}, we see that
$\otr_\tau B_*(s)$ has a simple pole at $s=1/2$ with residue given by
$$
-\frac12 \int_{\bX} \bar{u}_{k}^{\ast} \bar{v}_{k}^{\ast}\> dh_*,
$$
where $h_g = (\rho^2 g)|_\bX$ and similarly for $h_\tau$.  By Lemma 4.16 of
\cite{PP:2000} this is equal to $m_{1/2}^*$.  Thus the residue of $H(s)$ at $s=1/2$
is $m_{1/2}^g - m_{1/2}^\tau$, which gives the proper singularity for $D_{g,\tau}(s)$.
\end{proof}

The final step is the connection between $D_{g,\tau}(s)$ and $\det S_{g,\tau}(s)$.

\begin{lemma}
\label{detsr}
\begin{equation}
\det S_{g,\tau}(s)=e^{f_2(s)}\frac{D_{g,\tau}(1-s)}{D_{g,\tau}(s)}, \label{eq.rdet.prefunc}
\end{equation}
for a polynomial $f_2$ of degree at most two.
\end{lemma}

\begin{proof}
We will study the function
\[
F(s) = (2s-1)\otr_\tau [\hat R_g(s) - \hat R_g(1-s) - R_0(s) + R_0(1-s)].
\]
By \ref{eq.rdet.dlog} we see that
\[
\frac{d}{ds}\log[D_{g,\tau}(s)/D_{g,\tau}(1-s)] = F(s) + 2c(2s-1)
\]
We claim that
\begin{equation}\label{logdetsr}
\frac{d}{ds}\log\det(S_{g,\tau}(s)) = - F(s),
\end{equation}
which will complete the proof.  The formula (\ref{logdetsr}) is established by a rather long technical
calculation, which we defer to Appendix \ref{app.logdetsr}
\end{proof}

Combining Lemma \ref{sreldet} and Lemma \ref{detsr} gives the functional equation
\begin{equation}\label{eq.rdet.func}
\frac{D_{g,\tau}(1-s)}{D_{g,\tau}(s)} = e^{f(s)}\frac{P_g(1-s)}{P_g(s)}\frac{P_{\tau}(s)}{P_{\tau}(1-s)},
\end{equation}
where $f$ is a polynomial of degree at most four.  In conjunction with Lemma \ref{divdr}, this proves
Proposition \ref{prop.divisor}.

\subsection{Growth estimates on the relative determinant}

\label{sec.hadamard.estimates}

The results of the preceding subsection show that $W(s)=D_{g,\tau}(s)P_{\tau}(s)/P_g(s)$
is entire and zero-free. We now wish to estimate the growth of $\log W(s)$ as
$\left|  s\right|  \rightarrow\infty$ and show it to be polynomial. Using
(\ref{eq.h.def}) and Lemma \ref{lemma.dr.trace}, it is easy to see that
\begin{align}
\frac{d}{ds}\log W(s)  &  =H(s)+c(2s-1)\label{eq.h.dlog}\\
&  +P_{\tau}^{\prime}(s)/P_{\tau}(s)-P^{\prime}(s)/P_g(s).\nonumber
\end{align}
for some constant $c$.

It will suffice to show that the right-hand side of (\ref{eq.h.dlog}) has
polynomial growth as $\left|  s\right|  \rightarrow\infty$. By the estimate
(\ref{eq.r0.est}), for any $\delta>0$ there is a
countable collection of disjoint discs $\left\{  D_{j}\right\}  $ with the
properties that 
\begin{enumerate}
\item  $\mathcal{R}_g\cup\mathcal{R}_{\tau}\subset\cup_{j}D_{j}$
\item  $\operatorname*{dist}(s,\mathcal{R}_g
\cup\mathcal{R}_{\tau})\geq C\langle s\rangle^{-2-\delta}$ for every
$s\in\mathbb{C}\backslash(\mathcal{R}_g\cup\mathcal{R}_{\tau})$.
\end{enumerate} 
Here $\langle s\rangle=(1+|s|^{2})^{1/2}$. To prove that $W^{\prime
}(s)/W(s)$ has polynomial growth, it suffices by the maximum modulus theorem
to prove a polynomial growth estimate in $\mathbb{C}\backslash(\cup_{j}D_{j}
)$. Standard estimates on Hadamard products show that, on $\mathbb{C}
\backslash(\cup_{j}D_{j})$, the estimates
\[
\left|  P_{\tau}^{\prime}(s)/P_{\tau}(s)\right|  \leq C\langle s\rangle^{4+\delta}
\]
and
\[
\left|  P^{\prime}(s)/P_g(s)\right|  \leq C\langle s\rangle^{4+\delta}
\]
hold. Thus we need only show that the first right-hand term in
(\ref{eq.h.dlog}), which is $\frac{d}{ds}\log D_{g,\tau}(s)$ up to a linear
polynomial, grows polynomially in $\mathbb{C}\backslash(\cup_{j}D_{j})$. By
the functional equation (\ref{eq.rdet.func}) and these observations, it
suffices to show that, for some $\varepsilon>0$, $\frac{d}{ds}\log D_{g,\tau}(s)$
grows polynomially in the half-plane $\re(s)>\frac{1}{2}-\varepsilon$ with the
discs $D_{j}$ removed.

To this end, we first examine the behavior of
\begin{equation}
\left(\frac{1}{2s-1}\frac{d}{ds}\right)^{2}\log D_{g,\tau}(s)
=\tr \left[  \hat R_g(s)^{2}-R_\tau(s) ^{2}\right]  \label{eq.logD.d2}
\end{equation}
using the spectral shift representation.

\begin{lemma}
Let $\varepsilon>0$. For $\re(s) > \frac12 +\varepsilon$ and $|s|>1$ we can
estimate
\[
\left|  \tr \left[  \hat R_g(s)^{2}-R_\tau(s)^{2}\right]  \right|  \le
C_{\varepsilon}.
\]
\end{lemma}

\begin{proof}
By (\ref{zeta.realw}) and (\ref{xi.ss}) we can express the trace we are trying to estimate in
terms of the spectral shift function:
\[
\tr \left[  \hat R_g(s)^{2}-R_\tau(s)^{2}\right]  =-2\int_{\mu}^{\infty
}(\lambda+s(s-1))^{-3}\xi(\lambda)d\lambda.
\]
The integrand can be factored as
\[
(\lambda+s(s-1))^{-3}=\lambda^{-3}[1+s(s-1)/\lambda]^{-3}
\]
Since the assumptions on $s$ keep $s(s-1)$ bounded away from $(-\infty,-\mu]$,
we can estimate
\[
|1+s(s-1)/\lambda|^{-3}<C_{\varepsilon},
\]
uniformly for $\lambda\in\lbrack\mu,\infty)$ and $s$ as indicated. The result
then follows from (\ref{eq.xi.bd}). 
\end{proof}

Integration of (\ref{eq.logD.d2}) now yields the estimate
\[
\left|  \log D_{g,\tau}(s)\right|  \leq C\langle s\rangle
^{4}
\]
again for $\re(s)>\frac{1}{2}+\varepsilon$ and $\left|
s\right|  >1$. Let $S_{\varepsilon}$ be the open strip
$$
S_{\varepsilon} = \{1/2-\varepsilon<\re(s)<1/2+\varepsilon\}.
$$ 
Using the functional equation (\ref{eq.rdet.func}) and the estimates on Hadamard 
products, we can now bound
\[
\left|  \log H(s)\right|  \leq C\langle s\rangle
^{4+\varepsilon}
\]
for $S\in\mathbb{C}\backslash S_{\varepsilon}$.

To complete the proof we must estimate $|\log D_{g,\tau}(s)|$ in $S_{\varepsilon}$,
which is rather delicate because of the poles.  Fortunately, by the
Phragm\'{e}n-Lindel\"{o}f Theorem an exponential growth estimate in the strip
will suffice to extend the polynomial bounds.

\begin{lemma}
\label{lemma.strip} The estimate
\[
\left|  \log D_{g,\tau}(s)\right|  \leq C(\eta)\exp(\left|  s\right|  ^{2+\eta})
\]
holds for any $\eta>0$ and all $s$ with $s\in S_{\varepsilon}\backslash\left[
\cup_{j}D_{j}\right]  $.
\end{lemma}

The proof is quite technical and relies on a parametrix construction and estimates from 
Guillop\'e-Zworski \cite{GZ:1997}.  We will review this construction and prove Lemma \ref{lemma.strip} in
Appendix \ref{app.resolvent}.

Since $P_g(s)$ and $P_{\tau}(s)$ are of finite order,
and the Maximum Modulus Theorem can be used to fill in estimates in the disks
$D_{j}$, Lemma \ref{lemma.strip} gives us an estimate
\[
\log\left|  W(s)\right|  \leq C\exp(\left|  s\right|
^{2+\eta})
\]
for $s\in S_{\varepsilon}$. In $\mathbb{C}\backslash S_{\varepsilon}$, we have
\[
\left|  \log W(s)\right|  \leq C\langle s\rangle
^{4+\delta}.
\]
The Phragm\'{e}n-Lindel\"{o}f Theorem now applies to give us a polynomial
bound on $\log W(s)$ over the whole plane.  We have shown:

\begin{proposition}\label{prop.qpoly}
$W(s)=e^{q(s)}$ where $q(s)$ is a polynomial of degree at most 4.
\end{proposition}

\subsection{Asymptotics and order}

\label{sec.hadamard.asymptotics}

To show that $q(s)$ has order at most two, and that it depends only on the
eigenvalues and resonances of $\Delta_{g}$, let us rewrite the equation in
Theorem \ref{thm.det.r} as
\begin{equation}
\log D_{g,\tau}(s)=q(s)+\log P_g(s)-\log Z_{\tau}(s)+\log Z_{\infty}(s),
\label{eq.logdr}
\end{equation}
and consider the asymptotic behavior of each term as $\re(s)\rightarrow\infty$. By
Proposition \ref{prop.qpoly} $q(s)$ is polynomial of degree at most 4.  The
definition of $Z_{\tau}(s)$ in terms of the length spectrum (see (\ref{eq.zeta.euler})) 
together with the crude estimate that
\[
\#\{\gamma:l(\gamma)\leq t\}\leq Ce^{t}
\]
shows that
\[
|\log Z_{\tau}(s)|\leq Ce^{-c\re(s)}
\]
for $\re(s)>1$. This fact is rather crucial, since it shows that the
resonances of $\tau$ have no influence on the asymptotic expansion.

The expansion of $\log Z_{\infty}(s)$ has been given in (\ref{eq.zinf.exp}).
The asymptotics of the relative determinant can be obtained via the heat expansion.
This analysis is identical to the expansion of the zeta regularized
determinant in the compact case as given in \cite{Sarnak:1987}, with the
important exception that the first relative heat invariant $a_{1}=0$. For
$\re(s)\rightarrow\infty$,
\begin{equation}
\log D_{g,\tau}(s)\sim a_{0}s(s-1)[1-\log s(s-1)]+\sum_{j=2}^{\infty}
(j-2)!\>a_{j}z^{-j+1}. \label{eq.dr.exp}
\end{equation}

The following result will be critical to the consideration of isopolar
sequences of metrics, because (in contrast to the compact case) the relative
determinant $D_{g,\tau}(s)$ is \textit{not} a spectral invariant.

\begin{proposition}
\label{prop.spec.inv} The Euler characteristic $\chi(X)$, the polynomial
$q(s)$, and all of the relative heat invariants $a_{j}$ are determined by the
set of eigenvalues and resonances of $\Delta_{g}$.
\end{proposition}

\begin{proof}
The eigenvalues and resonances of $\Delta_{g}$ fix the Hadamard factor $P_g(s)$,
which must have an asymptotic expansion for $\re(s)\rightarrow\infty$ by the
analysis of the other terms in (\ref{eq.logdr}). Since $a_{1}=0$, the only
term of the form $\log s$ in the expansions is $\tfrac{1}{3}\chi(X)\log s$.
Therefore this term must be cancelled by a corresponding term in the expansion
of $\log P_g(s)$, from which we see that $\chi(X)$ is an isopolar invariant
(which implies $Z_{\infty}(s)$ is also). Knowing this, we observe that none of
the terms in the expansions of $q(s)$ and $\log D_{g,\tau}(s)$ could cancel with
each other and conclude that all coefficients are isopolar invariants.
\end{proof}

To complete the proof of Theorem \ref{thm.det.r} we note that $\log
P_g(s)=\mathcal{O}(|s|^{2}\ln|s|)$ as $s\rightarrow\infty$ with a similar
estimate for $P_{\tau}(s)$. On the other hand the heat expansion (\ref{eq.dr.exp}) 
can be used to derive an $\mathcal{O}(|s|^{2}\ln|s|)$ estimate for $\log
D_{g,\tau}(s) $ in a sector such as $|\arg(s)|\leq\varepsilon$. It follows that the
degree of $q(s)$ is actually two or less.

\section{Compactness for isopolar classes}

\label{sec.compact}

Let $\{g_k\}$ be a sequence of isopolar metrics on $X$ satisfying the hypotheses
of Theorem \ref{thm.compact}.  By pulling back by diffeomorphisms if necessary, 
we will assume
that the uniformizing metrics $\tau_k$ are arranged so that their convex cores
coincide.  That is, $\hat X_{\tau_k} = K$ for some fixed compact set $K\subset X$.
(We will see
below that such an alignment of the convex cores is necessary for convergence
of a subsequence.)
Then by assumption we have supp$(\varphi_k) \subset K$.
We wish to establish the existence of a convergent subsequence in the 
$\rho^{-2} \cinf(\bar X)$ topology.

The Hadamard factor depends only on the set of eigenvalues and resonances, so
we have a single function $P(s)$ for all $g_{k}$. Furthermore, Proposition
\ref{prop.spec.inv} tells us that the the polynomial $q(s)$ in the relation,
\begin{equation}
D_{g_k,\tau_k}(s)=e^{q(s)}\frac{P(s)}{Z_{\tau_k}(s)Z_{\infty}(s)}, \label{eq.dnr}
\end{equation}
is also independent of $k$.  (Recall that the
factor $Z_{\infty}(s)$ depends only on $\chi(X)$ and so has no $k$ dependence either.)
In addition to (\ref{eq.dnr}), we have also
the invariance of the relative heat invariants $a_{j} = a_j(\varphi_k,\tau_k)$.
And lastly, we have the Polyakov formula expressing $\log D_{g_k,\tau_k}(1)$ 
in terms of $\varphi_{k}$ and $\tau_k$  (but note that $\log D_{g_k,\tau_k}(1)$ is not
independent of $k$).

\begin{proof}[Proof of Theorem \ref{thm.compact}]
By combining the Polyakov formula (Proposition \ref{prop.polyakov}) with (\ref{eq.dnr}) we
obtain
\begin{equation}
-\frac{1}{6\pi}\int_{X} \left(\frac{1}{2}|\nabla_{\tau}\varphi_{k}|^{2}
-\varphi_{k}\right)\>d\tau_{k}=c-\log Z_{\tau_k}(1), \label{eq.poly.zn}
\end{equation}
where $c=q(1)+\log P(1)-\log Z_{\infty}(1)$ which is independent of $k$. 
Solving for $\varphi_k$ gives
$$
\int_X \varphi_k \>d\tau_k = \frac{1}{2}\int_{X} |\nabla_{\tau}\varphi_{k}|^{2}\>d\tau_{k}
+ 6\pi [c-\log Z_{\tau_k}(1)]
$$
Inspection of the definition (\ref{eq.zeta.euler}) of $Z_{\tau}$ yields that $- \log Z_{\tau_k}(1)\ge 0$,
and hence
\[
\int_{X}\varphi_{k}\>d\tau_{k}\geq 6\pi c.
\]
On the other hand, recall that the zeroth heat invariant,
\[
a_{0}=\frac{1}{4\pi}\int_{X}(e^{2\varphi_{k}}-1)\>d\tau_{k},
\]
is also a constant. An application of Jensen's inequality then gives
\[
\int_{X}\varphi_{k}\>d\tau_{k}\leq\frac{1}{2}\log(4\pi a_{0}+1).
\]
The left-hand side of (\ref{eq.poly.zn}) is thus bounded above and the
right-hand side is bounded below. The conclusion is a set of bounds:
\begin{equation}
6\pi c\leq\int_{X}\varphi_{\tau_{k}}\>d\tau_{k}\leq\frac{1}{2}\log(4\pi a_{0}+1),
\label{eq.phi.bnd}
\end{equation}
\begin{equation}
\int_{X}|\nabla_{\tau_k}\varphi_{k}|^{2}\>d\tau_{k}\leq\log(4\pi a_{0}+1)-12\pi c,
\label{eq.dphi.bnd}
\end{equation}
and
\begin{equation}
0\leq-\log Z_{\tau_k}(1)\leq\frac{1}{12\pi}\log(4\pi a_{0}+1)-c.
\label{eq.zn.bnd}
\end{equation}

By Corollary \ref{ZProperCoro}, the bound (\ref{eq.zn.bnd})
implies that a subsequence of the $\tau_{k}$'s converges, up to pull-back
by diffeomorphisms, in the $C^{\infty}$ topology on $X$ to some hyperbolic metric 
$\tau_{\infty}$.   Note that, by the construction made in the proof of Theorem \ref{ZProper}, 
we see that each of these diffeomorphisms preserves the respective pants decompositions.
Therefore the assumption $K= \hat X_{\tau_k}$ is preserved through the pull-backs.
Assume that a convergent subsequence has been chosen and diffeomorphisms applied to
give convergence in $\cinf(X)$.  For simplicity we will continue to denote the sequence
by $\tau_k$.
 
On the ends $F_j = (0,\infty)_t \times S^1_\theta$ these metrics can all be written as
$$
d\tau_k|_{F_j} = dt^2 + \ell_{k,j}^2 \cosh^2 t\>d\theta^2.
$$
The topology of $\rho^{-2} \cinf(\bar X)$ is independent of the choice of $\rho$.
With the choice $\rho = e^{-t}$, convergence of $\tau_k$ in $\rho^{-2} \cinf(\bar X)$ follows from 
the convergence $\ell_{k,j} \to\ell_{\infty,j}$ for all $j$.  

It remains to show that $\{\varphi_{k}\}$ has a convergent subsequence.
In the estimates above, the $C^{\infty}$ convergence of the 
$\tau_{k}$ allows us to replace $\tau_{k}$ by $\tau_{\infty}$ and still obtain uniform bounds.
From (\ref{eq.phi.bnd}) and (\ref{eq.dphi.bnd}) we derive estimates,
\[
\left|  \int_{K}\varphi_{k}\>d\tau_{\infty}\right|  \leq C,\qquad
\int_{K}|\nabla_{\tau_{\infty}}\varphi_{k}|^{2}\>d\tau_{\infty}\leq
C,
\]
which allow us to bound $\varphi_{k}$ uniformly in $H^{1}(X)$. From
this point on the analysis is exactly as in \S2 of Osgood-Phillips-Sarnak
\cite{OPS:1988b}. The constancy of the second relative heat invariant
(\ref{eq.a2}), in conjunction with the estimates above, is used to bound
$\varphi_{k}$ uniformly in $C^{0}(X)$. Then the higher heat invariants
(\ref{eq.aj}) allow a bootstrap argument extending the uniform bounds to
$H^{m}(X)$ for all $m$. This, together with the restriction $\varphi_k \in \cinf_0(K)$, 
implies compactness of the $\varphi_k$'s in the $C^{\infty }_0(K)$ topology, and
completes the proof.
\end{proof}

\appendix

\section{Finiteness and properness for hyperbolic surfaces}
\label{sec.proper}

Let $S$ be an oriented compact topological surface $S$ with boundary $\partial S$.
Let ${\mathcal M}(S)$ denote the set of equivalence classes (up to isometry)
hyperbolic metrics $h$ on $S$ such that $\partial S$ is geodesic with respect 
to $h$. 

We will use the $C^{\infty}$ topology on ${\mathcal M}(S)$:
A sequence $\{h_n\} \subset {\mathcal M}(S)$ converges if 
and only if there exists $h'_n \in h_n$ such that the coordinate components 
of $h'_n$ converge in $C^{\infty}$.

The (primitive) length spectrum $\Lambda(h)$ 
does not depend on the choice of representative for the class $h$.
Hence the class $h$ defines through $(\ref{eq.zeta.euler})$ a dynamical zeta
function $Z_h(s)$.
Note also that the total length, $\ell_{h}(\partial S)$, of the boundary
does not depend on the representative of $h$.

\begin{theorem} \label{ZProper}
Let $R>0$. The set of all $h \in {\mathcal M}(S)$ 
such that $-\log Z_{h}(1) \leq R$ and $\ell_{h}(\partial S) \leq R$
is compact.
\end{theorem}

The geometry of a convex co-compact hyperbolic surface $(X,\tau)$  
is determined by the geometry of its convex core. Moreover,
the zeta function associated to $(X,\tau)$ equals the 
zeta function  associated to its convex core. Therefore, we have the
following:

\begin{corollary} \label{ZProperCoro}
Let $X$ be an open surface of finite topological type,
with $\tau_n$ a sequence of complete hyperbolic metrics
on $X$.
If $-\log(Z_{\tau_n}(1))$ is bounded from above, then there
exists a hyperbolic metric $\tau_\infty$ and a sequence of 
diffeomorphisms $\phi_n: X \rightarrow X$ such that
$\phi_n^*(\tau_n)$ converges to $\tau_\infty$ in the $C^{\infty}$ 
topology.
\end{corollary}

Before considering the proof of Theorem \ref{ZProper}, 
we give provide an example that shows that an upper bound on  
$\ell(\partial S)$ is a necessary condition.

\begin{example} \label{PantsExample}
Let $S$ be the sphere with three open discs removed, a topological
{\em pair of pants}.
Then each $h \in {\mathcal M}(S)$ is determined by the
lengths $\ell(\gamma_1),\ell(\gamma_2), \ell(\gamma_3))$,
of the connected components  $\gamma_1,\gamma_2,\gamma_3$
of  $\partial S$. Conversely, given 
$(\ell_1, \ell_2, \ell_3) \in ({\mathbb R}^+)^3$ there
exists $h(\ell_1,\ell_2, \ell_3)$ such that $\ell(\gamma_i)= \ell_i$.
(See, for example, Theorem 3.1.7 \cite{Buser:1992}).

We claim that 
\begin{equation} \label{Lower}
   \inf \Lambda(h(\ell_1,\ell_2, \ell_3)) = \inf\{ \ell_1,\ell_2, \ell_3\} 
\end{equation}
and hence if $\ell_1,\ell_2, \ell_3$ are bounded from below then 
$-\log(Z_{h(\ell_1,\ell_2, \ell_3)})$ is bounded from above.
On the other hand, the family $h(\ell,\ell, \ell)$, for example, has no limit
as $\ell$ tends to infinity.  Hence, it is necessary to assume
an upper bound on $\ell(\partial S)$ in Theorem \ref{ZProper}. 

To verify (\ref{Lower}) we note that any non-null homotopic 
simple closed curve $\alpha$ on $S$ 
is homotopic to either $\gamma_1$, $\gamma_2$, or $\gamma_3$. Indeed,
by the Jordan curve theorem, $S\setminus \alpha$ 
has two components, and it follows that $\alpha$ is homotopic to one of the 
$\gamma_i$. Any closed geodesic $\beta$ 
has finitely many self-intersections and they are all transverse. 
By doing a surgery at each crossing, one obtains a finite number 
of simple closed curves each of whose length is greater 
than $\inf\{ \ell_1,\ell_2, \ell_3\}$. Thus, 
$\ell(\beta) \geq \inf\{ \ell_1,\ell_2, \ell_3\}$ and the
claim follows.
\end{example}

The proof of Theorem \ref{ZProper} relies on the following
variant of Bers' Theorem.

\begin{theorem}[Bers' Theorem for surfaces with geodesic boundary] \label{BersThm}
Let $S$ be a compact surface with boundary. There exists a constant $c=c(S)$
such that for any hyperbolic metric $h$ on $S$ inducing geodesic boundary,
there exists a decomposition of $S$ into pairs of pants $P_1, \ldots, P_n$ 
such that $\ell_h(\partial P_i) \leq c$ for each $i$.
\end{theorem}

\begin{proof}
The idea of the proof comes from Theorem 5.2.3 of \cite{Buser:1992}. 
Let $\exp: [0, \infty] \times \partial S \rightarrow S$ be the 
exponential map associated to the normal bundle of $\partial S$.
Let $t$ be the largest $t$ such that the restriction of 
$\exp$ to $[0, t[ \times \partial S$ is injective.
Note that on each connected component of $K_i$ of 
$[0, \infty] \times \partial S$ we have 
$\exp^*(h)= d \rho ^2 + \ell_i^2 \cosh(\rho)^2 d \theta_i^2$ 
where $\ell_i$ is the length of $\partial K_i$. 
Hence, 
\begin{eqnarray*}
    \ell(\exp(\{t\} \times \partial S)) &=& \cosh(t) \cdot \ell(\partial S) \\
   \Area(\exp([0, t] \times \partial S)) &=& \sinh(t) \cdot \ell(\partial S).    
\end{eqnarray*}
From $\sinh^2(t) +1 = \cosh^2(t)$, we then find that 
\begin{equation} \label{Curve}
 \ell^2(\exp(\{t\} \times \partial S))~ \leq~ \ell^2(\partial S) + (\Area(S))^2.
\end{equation}

The set $\exp(\{t\} \times \partial S)$ is a union of simple closed curves
each of which is freely homotopic to a unique simple closed geodesic $\gamma_i$.
Let $P_1, \ldots P_k$ be the connected components of $S \setminus (\cup \gamma_i)$
that are pairs of pants.  By construction, $k \geq 1$ unless $S$ is a pair
of pants---see p. 126-129 \cite{Buser:1992}---and  $\ell(S \setminus (\cup P_i))
   \leq \ell(\cup \gamma_i) \leq \ell(\exp(\{t\} \times \partial S))$.
Therefore, by (\ref{Curve})
\begin{equation*}
   \ell^2(\partial (S \setminus (\cup P_i)))~ \leq~ 
   \ell^2(\partial S) + (\Area(S))^2.
\end{equation*}
Since $k \geq 1$, we have $\Area(S \setminus (\cup P_i)) \leq \Area(S) - 2 \pi$
and since the claim is vacuous
for a pair of pants, the general claim follows by induction.
\end{proof}

\begin{proof}[Proof of Theorem \ref{ZProper}]
Suppose that $-\log(Z_{h}(1)) \leq R$.
Since each term of (\ref{eq.zeta.euler}) is positive and decreasing in $s>0$, 
we find that $-\log(1-\exp(-(2+k)\cdot \ell(\gamma))) <R$ 
for all $\gamma$ and $k$.  In particular, for each 
primitive---and hence every---closed geodesic $\gamma$
\begin{equation} \label{MumLower}
     \ell(\gamma) \geq \varepsilon_R >0
\end{equation}
where $\varepsilon_R = -\frac{1}{2} \ln(1-\exp(-R))$. 

Let $h_n$ be a sequence of hyperbolic metrics on $S$ inducing 
geodesic boundary and satisfying $\inf \Lambda(h_n) \geq \varepsilon_R$
and $\ell(\partial S)<R$. It suffices to show that there exists
a subsequence, still called $h_n$, and 
diffeomorphisms $\psi_n: S \rightarrow S$ such that
$\psi_n^*(h_n)$ converges in $C^{\infty}$.

Applying Theorem \ref{BersThm} to each 
$h_n$ gives an infinite sequence of pants decompositions 
$\{P_1^n, \ldots, P_k^n\}$ with $\ell_{h_n}(\partial P_i) <c$
for all $n$ and $i$.
Since there are only finitely many combinatorial types of pants
decompositions---see \cite{Buser:1992} \S 3.6---we may assume without 
loss of generality that $\{P_1^n, \ldots, P_k^n\}$ has constant 
combinatorial type. It follows that there exist diffeomorphisms 
$\phi_n: S \rightarrow S$ such that $\phi_n(P_i^n)=P_i^1$ for
$i=1, \ldots k$. Hence, by pulling back $h_n$ by $\phi_n$,
we may assume without loss of generality that each $h_n$
gives the same pants decomposition $\{P_1, \ldots, P_n\}$.  

By Theorem  \ref{BersThm} and (\ref{MumLower}), 
we have $\varepsilon_R \leq \ell_{h_n}(\gamma_{ij}) \leq c$
for each boundary component $\gamma_{ij}$ of $\partial P_i$,
$j=1,2,3$. Therefore there is a subsequence, still denoted $h_n$, 
such that each $\ell_{h_n}(\gamma_{ij})$ converges as $n$ tends
to infinity.  Thus, by the discussion in Example \ref{PantsExample}
and the reference given there, there exist diffeomorphisms 
$\psi_n^i: P_i \rightarrow P_i$ such that 
$\psi_n^i{*}(h_n|_{P_i})$ converges to a metric on $P_i$ in $C^{\infty}$.

Finally, by perturbing each $\psi_n^i$ in a collar neighborhood of
each boundary component of $P_i$, one can construct a suitable
diffeomorphism of the entire
surface $\psi_n: S \rightarrow S$.  Moreover, the associated {\em twist
angles} (see \cite{Buser:1992} \S 3.3.)
are bounded, and hence by passing to a further subsequence if necessary,
we may assume that $\psi_n^*(h_n)-(\psi_n^i)^*(h_n)$ is Cauchy.
It follows that  $\psi_n^*(h_n)$ converges in $C^{\infty}$ to a metric $h$ on $S$.  
\end{proof}

Our variant of Bers' theorem also has the following corollary.
Compare with Theorem 13.1.3 in \cite{Buser:1992}.

\begin{theorem} \label{ZFinite}
Let $h_0 \in {\mathcal M}(s)$ and $R>0$.
Then the set of all $h \in {\mathcal M}(S)$ 
such that $\Lambda(h)= \Lambda(h_0)$ and $\ell_{h}(\partial S) \leq  R$ 
is finite.
\end{theorem}

\begin{proof}
The set in question is compact by Theorem \ref{ZProper}.
The twist angles are determined by $\Lambda(h_n)$
as are the lengths of the boundary components of the $P_i$.
Hence the set is discrete as well as compact.
\end{proof}

\section{Resolvent construction and estimates}

\label{app.resolvent}

We briefly review the construction of the resolvent carried out in Guillop\'e-Zworski
\cite{GZ:1997} since we will used some detailed information from that
construction to prove key estimates on the relative determinant. We follow
closely the outline of \cite{GZ:1997}.
The results of this section apply to any metric that is hyperbolic near infinity,
but not to a general asymptotically hyperbolic metric.
For notational simplicity, let us assume $\tau$ is a hyperbolic metric on
$X$ and $g$ is a perturbation that is equal to $\tau$ on $\hat X$.

The cylinders $F_{j}$ with the hyperbolic metric $\tau$ are isometric to
hyperbolic half-cylinders $F_{j}^{0}=(0,\infty)_{t}\times S_{\theta}^{1}$
with metric
\begin{equation}
\tau_{j}^{0}=dt^{2}+\ell_{j}^{2}\cosh^{2}t\,\,d\theta^{2}
\label{eq.funnel.metric}
\end{equation}
where $\ell_{j}$ is the geodesic length of the circle at $t=0$ (which is a
closed geodesic). If $\Delta_{F_{j}^{0}}$ is the Laplacian on $(F_{j}
^{0},\tau_{j}^{0})$, we let
\[
R_{F_{j}^{0}}(s)=(\Delta_{F_{j}^{0}}-s(1-s))^{-1}
\]
be the resolvent for the hyperbolic Laplacian on $F_{j}^{0}$ with Dirichlet
boundary conditions at $t=0$. This resolvent can be computed explicitly
(see for example \cite{Guillope:1990} or \cite{GZ:1995a}) and is known to have
poles contained in the set of $\zeta_{n,k}=-k+2\pi in/\ell_{j}$ where $n$ is
any integer and $k=0,1,2,\cdots$. In particular, $R_{F_{j}^{0}}(s)$ is entire
in any half-plane $\re(s)>\varepsilon$. Finally, it follows from explicit
formulas that if $\chi$ and $\psi$ are smooth, compactly supported functions
in $F_{j}^{0}$ with disjoint supports, $\chi R_{F_{j}^{0}}(s)\psi$ has a
smooth kernel with derivatives bounded uniformly in $s$ with $\re
(s)>\varepsilon$ for any fixed $\varepsilon>0$.

First we describe the parametrix constructed from model operators. Let
$\eta\in\cinf(\mathbb{R})$ with $\eta(t)=1$ for $t<1/3$,
$\eta(t)=0$ for $t>2/3$, and let $\eta_{a}(t)=\eta(t-a)$. We will pick $a>1$
in what follows. Let $\chi_{a}\in\cinf_0(X)$ with
$\chi_{a}=1$ on $Z$ and $\chi_{a}=\eta_{a}$ on each funnel $F_{j}$
(referring to the coordinates $(t,\theta)$ as in (\ref{eq.funnel.metric})). We
denote by $\chi_{a,j}$ the restriction of $\chi_{a}$ to $F_{j}$. For a
fixed $a>0$ and real $s_{0}$ with $s_{0}>1$, we set
\begin{equation}
E_{g}(s)=Q_{0}(s_{0})+Q(s) \label{eq.e.def}
\end{equation}
where
\begin{equation}
Q_{0}(s_{0})=\chi_{a+2}R_{g}(s_{0})\chi_{a+1} \label{eq.q0.def}
\end{equation}
and
\begin{equation}
Q(s)=\sum_{j=1}^{M}(1-\chi_{a,j})J_{j}^{\ast}R_{F_{j}^{0}}(s)J_{j}
(1-\chi_{a+1,j}). \label{eq.q.def}
\end{equation}
Here $J_{j}:F_{j}\rightarrow F_{j}^{0}$ is an isometry mapping the funnel end
to the model manifold.

It is not difficult to compute that
\[
(\Delta_{g}-s(1-s))E_{g}(s)=I+L(s_{0},s)
\]
where
\begin{align}
L(s_{0},s)  &  =(s_{0}(1-s_{0})-s(1-s))Q_{0}(s_{0})+\left[  \Delta_{g}
,\chi_{a+2}\right]  R_{g}(s_{0})\chi_{a}\label{eq.l.def}\\
&  -\sum_{j=1}^{M}\left[  \Delta_{g},\chi_{a}\right]  J_{j}^{\ast}
R_{F_{j}^{0}}(s)J_{j}(1-\chi_{a+1,j})\nonumber
\end{align}
The first right-hand term is a Hilbert-Schmidt integral operator with
compactly supported kernel, while the second is a compactly supported operator
with smooth kernel owing to the fact that the derivatives of $\chi_{a+2}$
and $\chi_{a}$ have disjoint supports. The third term in $L(s_{0},s)$ is a
sum of operators with kernels belonging to $(\rho^{\prime})^{s}\cinf
(\bar{X}\times\bar{X})$ having compact
support in the first variable (here $\rho^{\prime}$ is a defining function for 
$\partial\bar{X}$ in the
second variable). To construct the resolvent $R_g(s)$ in the half plane
$\re(s)>1/2-N$, one inverts the operator $(I+L(s_{0},s))$ viewed as a map from
$\rho^{N}L^{2}(X)$ to itself, using the analytic Fredholm theorem.

One then obtains
\begin{align}
R_g(s)  &  =E_{g}(s)(I+L(s_{0},s))^{-1}\label{eq.r}\\
&  =E_{g}(s)-E_{g}(s)(I+L(s_{0},s))^{-1}L(s_{0},s).\nonumber
\end{align}

\begin{proof}[Proof of Lemma \ref{lemma.strip}]
It suffices to show that a similar estimate holds for $H(s)$ as defined in
(\ref{eq.h.def}) since the desired estimate follows by integration.  From
(\ref{eq.r}) we may write
\begin{equation}\label{eq.rdiff}\begin{split}
R_g(s)-R_\tau(s)  &  =E_{g}(s)-E_{\tau}(s)\\
& \qquad +(E_{g}(s)-E_{\tau}(s))(I+L(s_{0},s))^{-1}L(s_{0},s)\\
& \qquad +E_{\tau}(s)(I+L(s_{0},s))^{-1}(L(s_{0},s)-L_{0}(s_{0},s))\\
& \qquad +E_{\tau}(s)\left[  (I+L(s_{0},s))^{-1}-(I+L_{0}(s_{0},s))^{-1}\right]
L_{0}(s_{0},s)\\
&  =T_{1}(s)+T_{2}(s)+T_{3}(s)+T_{4}(s).\nonumber
\end{split}\end{equation}
\newline It suffices to estimate the zero-trace of the quantities
$T_{i}(s)-T_{i}(s_{0})$. As we will see, the difference is important only when
$i=1,2$ and otherwise it is possible to prove polynomial bounds on the
zero-traces of $T_{i}(s)$, $i=3,4$. Indeed, $T_{1}(s)-T_{1}(s_{0})$,
$T_{2}(s)-T_{2}(s_{0})$, and $T_{3}(s)$ are operators with continuous kernel
whose support is compact in at least one of the variables; thus the zero-trace
is actually a trace:
\begin{align*}
\otr_\tau \left[  T_{i}(s)-T_{i}(s_{0})\right]   &  =\tr \left[
\chi_{a+3}\left[  T_{i}(s)-T_{i}(s_{0})\right]  \chi_{a+3}\right]
,\,\,i=1,2\\
\otr_\tau \left[  T_{3}(s)\right]   &  =\tr \left[  \chi_{a+3}T_{3}(s)\chi_{a+3}\right].
\end{align*}
Thus it suffices to estimate the trace norms of the operators
\[
\chi_{a+3}\left[  T_{1}(s)-T_{1}(s_{0})\right]  \chi_{a+3},\quad
\chi_{a+3}\left[  T_{2}(s)-T_{2}(s_{0})\right]  \chi_{a+3},\quad
\chi_{a+3}T_{3}(s)\chi_{a+3}.
\]
We will use a separate argument to estimate $T_{4}$. In what follows,
$\left\|  \,\cdot\,\right\|  $, $\left\|  \,\cdot\,\right\|  _{1}$, and
$\left\|  \,\cdot\,\right\|  _{2}$ will denote respectively the operator norm,
trace norm, and Hilbert-Schmidt norm.

The $T_{1}(s)$ term
is, by (\ref{eq.e.def}), (\ref{eq.q0.def}), and (\ref{eq.q.def}), equal to
\[
\chi_{a+2}\left[  R_g(s_{0})-R_\tau(s_{0})\right]  \chi_{a+1}
\]
since the model operators occuring in $Q(s)$ are the same for the two
problems. It follows that $T_{1}(s)-T_{1}(s_{0})=0$ so this term makes no
contribution.

To estimate the remaining terms, we first note
that
\[
(I+L(s_{0},s))^{-1}L(s_{0},s)\chi_{a+3}=(I+K(s_{0},s))^{-1}K(s_{0},s)
\]
where
\[
K(s_{0},s)=L(s_{0},s)\chi_{a+3}
\]
and similarly for $L$ and $K$ replaced by $L_{0}$ and $K_\tau $. It is not
difficult to see that
\begin{equation}
\left\|  K(s_{0},s)\right\|  _{2}\leq C_{\varepsilon}\langle s\rangle^{2}
\label{eq.k.est}
\end{equation}
for $s\in S_{\varepsilon}$, with a similar estimate for $K_\tau $. From
\cite{GZ:1997}, Lemma 3.6, we have the estimate
\begin{equation}
\left\|  (I+K(s_{0},s))^{-1}\right\|  \leq\exp(C_{\eta}\langle s\rangle
^{2+\eta}) \label{eq.kfred.est}
\end{equation}
for any $\eta>0$ and all $s\in S_{\varepsilon}\backslash\cup_{j}D_{j}$ (The
statement of \cite{GZ:1997}, Lemma 3.6 also excludes singularities of the
model resolvents on the half-cylinders $F_{j}^{0}$, but these singularities
lie in the half-plane $\re(s)\leq0$ and so are already excluded from
$S_{\varepsilon}$).

We now estimate the remaining terms. First of all,
\begin{align*}
\chi_{a+3}T_{2}(s)\chi_{a+3}  &  =\chi_{a+2}\left[  R_g(s_{0}
)-R_\tau(s_{0})\right]  \chi_{a+1}\times\\
&  (I+K(s_{0},s))^{-1}K(s_{0},s)
\end{align*}
but the difficulty here is that the difference
\[
R_g(s_{0})-R_\tau(s_{0})
\]
is not trace-class. On the other hand, the difference
\begin{align*}
\chi_{a+3}\left[  T_{2}(s)-T_{2}(s_{0})\right]  \chi_{a+3}  &
=\chi_{a+2}\left[  R_g(s_{0})-R_\tau(s_{0})\right]  \chi_{a+1}\times\\
&  (I+K(s_{0},s_{0}))^{-1}\left[  K(s_{0},s)-K(s_{0},s_{0})\right]
(I+K(s_{0},s))^{-1}
\end{align*}
and the factor $K(s_{0},s)-K(s_{0},s_{0})$ is easily seen to be a trace-class
operator whose trace norm has at most polynomial growth in $s$.  It now follows
from (\ref{eq.k.est}), (\ref{eq.kfred.est}), and the fact that the first
factor is a fixed trace-class operator that
\[
\left\|  \chi_{a+3}\left[  T_{2}(s)-T_{2}(s_{0})\right]  \chi_{a+3}
\right\|  \leq\exp(C_{\eta}^{\prime}\langle s\rangle^{2+\eta}).
\]
To bound $T_{3}$, we note that
\[
\chi_{a+3}T_{3}(s)\chi_{a+3}=(\chi_{a+3}E_{\tau}(s)\chi_{a+4}
)(I+K(s_{0},s))^{-1}\left[  K(s_{0},s)-K_\tau (s_{0},s)\right]
\]
and use the fact that
\[
\left\|  ABC\right\|  _{1}\leq\left\|  A\right\|  _{2}\left\|  B\right\|
\left\|  C\right\|  _{2}
\]
together with the bounds
\begin{equation}
\left\|  \chi_{a+3}E_{\tau}(s)\chi_{a+4}\right\|  _{2}\leq C\langle
s\rangle^{2}, \label{eq.e0.est}
\end{equation}
(\ref{eq.k.est}), and (\ref{eq.kfred.est}) to conclude that
\[
\left\|  \chi_{a+3}T_{3}(s)\chi_{a+3}\right\|  _{1}\leq\exp(C_{\eta
}^{\prime}\langle s\rangle^{2+\eta}).
\]

To bound $T_{4}(s)$ we make a slightly different argument since the integral
kernel need not have compact support. Instead, we consider $\chi_{\varepsilon
}T_{4}(s)\chi_{\varepsilon}$ where $\chi_{\varepsilon}$ is the characteristic
function of the set $\rho\geq\varepsilon$. Write $E_{\tau}$ for $E_{\tau}(s)$, $L$
for $L(s_{0},s)$ and similarly $L_{0}$ for $L_{0}(s_{0},s))$. We compute, by
cyclicity of the trace,
\begin{align}
\tr (\chi_{\varepsilon}T_{4}(s)\chi_{\varepsilon})  &
=\tr (\chi_{\varepsilon}E_{\tau}\left[  (I+L)^{-1}-(I+L_{0}
)^{-1}\right]  L_{0}\chi_{\varepsilon})\label{eq.t4}\\
&  =\tr ((I+L_{0})^{-1}L_{0}\chi_{\varepsilon}^{2}
E_{\tau}(I+L)^{-1}\left[  L_{0}-L\right]  )\nonumber\\
&  =\tr (\chi_{a+5/2}L_{0}(I+L_{0})^{-1}L_{0}\chi
_{\varepsilon}^{2}E_{\tau}(I+L)^{-1}L_{0}\chi_{a+3})\nonumber\\
& \qquad \tilde{h} -\tr (\chi_{a+5/2}(I+L_{0})^{-1}L_{0}\chi_{\varepsilon
}^{2}E_{\tau}(I+L)^{-1}L\chi_{a+3})\nonumber
\end{align}
We have inserted the cutoff functions $\chi_{a+5/2}$ on the left owing to
the mapping properties of $L$ and $L_{0}$ and on the right by cyclicity of the
trace. Observe that each of the terms in the last two lines of (\ref{eq.t4})
is actually trace-class since $L\psi$ and $L_{0}\psi$ are Hilbert-Schmidt for
any function $\psi\in\cinf_0(X)$. Moreover, $L_{0}
\chi_{\varepsilon}^{2}E_{\tau}=L_{0}\chi_{a+3}E_{\tau}$ for sufficiently small
$\varepsilon$ since $\chi_{a+3}E_{\tau}=\allowbreak E_{\tau}$, so that
\begin{align*}
\otr_\tau (T_{4}(s))  &  =\tr (\chi_{a+5/2}L_{0}(I+L_{0}
)^{-1}L_{0}\chi_{a+3}E_{\tau}(I+L)^{-1}L_{0}\chi_{a+3})\\
& \qquad +\tr (\chi_{a+5/2}(I+L_{0})^{-1}L_{0}\chi_{a+3}
E_{\tau}(I+L)^{-1}L\chi_{a+3})\\
&  =\tr (\chi_{a+3}L_{0}(I+K_\tau )^{-1}K_\tau \chi_{a+4}
E_{\tau}(I+K)^{-1}K_\tau )\\
& \qquad +\tr (\chi_{a+3}(I+K_\tau )^{-1}K_\tau \chi_{a+4}
E_{\tau}(I+K)^{-1}K).
\end{align*}
We can now use (\ref{eq.k.est}) and (\ref{eq.kfred.est}) together with
(\ref{eq.e0.est}) to conclude that
\[
\left|  \otr_\tau (T_{4}(s))\right|  \leq\exp(C_{\eta}^{\prime}\langle
s\rangle^{2+\eta})
\]
and complete the proof.
\end{proof}

\section{Logarithmic derivative of the relative scattering operator}

\label{app.logdetsr}

To prove the formula (\ref{logdetsr}) we need to explicitly evaluate the
0-traces appearing in $F(s)$. We will employ the basic strategy used in \S6.1
of \cite{PP:2000}; similar calculations are done in Proposition 4.1 of
\cite{GZ:1997}. For the hyperbolic metric $\tau$ we have an identity:
\begin{equation}
R_\tau(s;z,z^{\prime})-R_\tau(1-s;z,z^{\prime})=\frac{1}{2s-1}
\int_{\partial X}E_{\tau}(1-s;z,y)E_{\tau}(s;z^{\prime},y)dh(y), \label{eq.rree}
\end{equation}
where $dh$ is the Riemannian density induced on ${\partial X}$ by $\rho
^{2}\tau|_{\partial X}$.  Note that the same density is induced by $\rho
^{2}g|_{\partial X}$.
Here $E_{\tau}(s;z,y)$ is the Poisson kernel for
$\Delta_\tau+s(s-1)$, which can be realized as a limit
\[
E_{\tau}(s;z,y)=(2s-1)2^{2s-1}\frac{\Gamma(s-1/2)}{\Gamma(1/2-s)}\rho
(z^{\prime})^{-s}R_\tau(s;z,z^{\prime})|_{z^{\prime}=y}.
\]
If $\hat R_g(s;z,z)$ is the kernel of the resolvent of $\hat\Delta_g$, with respect 
to the measure $d\tau$, and $\hat E_{g}(s;z,y)$ is defined by a limit as above, then it 
is simple to check that the corresponding relation holds.

With these relations we can write
\begin{equation}
\begin{split}
F(s)  &  ={\FP_{\varepsilon\downarrow0}}\int_{\rho\geq\varepsilon}
\int_{\partial X}\Bigl[\hat E_{g}(1-s;z,y) \hat E_{g}(s;z,y)\\
&  \hskip1in-E_{\tau}(1-s;z,y)E_{\tau}(s;z,y)\Bigr]dh(y)d\tau(z).
\end{split}
\label{eq.Fintegral}
\end{equation}
Technically we should assume here that $\rho$ is a suitable defining function 
for the definition of $\otr_\tau$.  However, the calculation will show that
$F(s)$ does not actually depend on the choice of defining function.

The next step is to transform the integral using the Maass-Selberg relation,
which for $E_{\tau}(s)$ reads:
\[
\begin{split}
&  (2s-1)\int_{\rho\ge\varepsilon} \left[  \int_{\partial X} E_{\tau}(1-s; z,y)
E_{\tau}(s; z,y) dh(y) \right]  d\tau(z)\\
&  \qquad= - \int_{\rho=\varepsilon} \int_{\partial X} \Bigl[ E_{\tau}(1-s; z,y)
\partial_{s} \partial_{\nu}E_{\tau}(s; z,y)\\
&  \hskip2in- \partial_{\nu}E_{\tau}(1-s; z,y) \partial_{s} E_{\tau}(s; z,y) \Bigr]
dh(y) d\sigma_{\varepsilon}(z),
\end{split}
\]
where $d\sigma_{\varepsilon}$ is the measure induced on $\{\rho=\varepsilon\}$ by $\tau$,
and $\partial_{\nu}$ is the inward Riemannian normal to $\{\rho=\varepsilon\}$.

Since the corresponding equation for $\hat\Delta_g$ contains extra factors, we will go
through the proof for this case. Let $\omega(s,t)=(s+t)(1-s-t)-s(1-s)$. The
first step is to use the property $(\hat\Delta_g-s(1-s))\hat E_{g}(s;\cdot,y)=0$ to write
\[
\begin{split}
&  (2s-1)\int_{\rho\geq\varepsilon}\left[  \int_{\partial X}
\hat E_{g}(1-s;z,y)\hat E_{g}(s;z,y)dh(y)\right]  d\tau(z)\\
&  \qquad=\lim_{t\rightarrow0^{+}}\frac{2s-1}{\omega(s,t)}\int_{\rho
\geq\varepsilon}\int_{\partial X}\Bigl[\hat E_{g}(1-s;z,y)\hat\Delta_g \hat E_{g}(s+t;z,y)\\
&  \hskip2in-\hat\Delta_g \hat E_{g}(1-s;z,y)\hat E_{g}(s+t;z,y)\Bigr]dh(y)d\tau(z).
\end{split}
\]
Then integrate by parts, using $\hat\Delta_g=e^{-\varphi}\Delta_\tau e^{-\varphi}$:
\[
\begin{split}
&  \qquad=\lim_{t\rightarrow0^{+}}\frac{2s-1}{\omega(s,t)}\int_{\rho=\varepsilon
}\int_{\partial X}\Bigl[e^{-\varphi(z)}\hat E_{g}(1-s;z,y)\partial_{\nu}(e^{-\varphi
(z)}\hat E_{g}(s+t;z,y))\\
&  \hskip1in-\partial_{\nu}(e^{-\varphi(z)}\hat E_{g}(1-s;z,y))e^{-\varphi
(z)}\hat E_{g}(s+t;z,y)\Bigr]dh(y)d\sigma_{\varepsilon}(z).\\
&  \qquad=-\int_{\rho=\varepsilon}\int_{\partial X}\Bigl[e^{-\varphi
(z)}\hat E_{g}(1-s;z,y)\partial_{\nu}(e^{-\varphi(z)}\partial_{s}\hat E_{g}(s;z,y))\\
&  \hskip1in-\partial_{\nu}(e^{-\varphi(z)}\hat E_{g}(1-s;z,y))e^{-\varphi(z)}
\partial_{s}\hat E_{g}(s;z,y)\Bigr]dh(y)d\sigma_{\varepsilon}(z).
\end{split}
\]

Since we know that $(\partial_{\nu}-s)E_{\ast}(s;z,y)$ is of lower order as
$z\rightarrow{\partial X}$, it is useful to extract terms of this form as we
substitute the Maass-Selberg relations back into (\ref{eq.Fintegral}). For
this purpose we note that $\partial_{s}\partial_{\nu}=\partial_{s}
(\partial_{\nu}-s)+s\partial_{s}+1$. We will break the resulting formula for
$F(s)$ up into pieces:
\[
F(s)=-\frac{1}{2s-1}\int_{\partial X} \left({\FP_{\varepsilon\downarrow0}}
\int_{\rho=\varepsilon}[J_{1}(z,y)+J_{2}(z,y)+J_{3}(z,y)]d\sigma_{\varepsilon}(z) \right)
dh(y),
\]
where
\[
\begin{split}
J_{1}(s;z,y)  &  =e^{-2\varphi(z)}\hat E_{g}(1-s;z,y)\hat E_{g}(s;z,y)-E_{\tau}(1-s;z,y)E_{\tau}
(s;z,y)\\
J_{2}(s;z,y)  &  =e^{-\varphi(z)}\hat E_{g}(1-s;z,y)\partial_{s}(\partial_{\nu
}-s)e^{-\varphi(z)}\hat E_{g}(s;z,y)\\
&  \qquad-(\partial_{\nu}-1+s)e^{-\varphi(z)}\hat E_{g}(1-s;z,y)e^{-\varphi(z)}
\partial_{s}\hat E_{g}(s;z,y)\\
&  \qquad-E_{\tau}(1-s;z,y)\partial_{s}(\partial_{\nu}-s)E_{\tau}(s;z,y)\\
&  \qquad+(\partial_{\nu}-1+s)E_{\tau}(1-s;z,y)\partial_{s}E_{\tau}(s;z,y)\\
J_{3}(s;z,y)  &  =(2s-1)e^{-2\varphi(z)}\hat E_{g}(1-s;z,y)\partial_{s}\hat E_{g}(s;z,y)\\
&  \qquad-(2s-1)E_{\tau}(1-s;z,y)\partial_{s}E_{\tau}(s;z,y)\\
\end{split}
\]
Let $I_{j}(s)$ denote the contribution to $F(s)$ from $J_{j}$. We will show
that $I_{1}(s)=I_{2}(s)=0$, while $I_{3}(s)$ is equal to the stated result.

To handle $J_{1}$ we reverse the identity (\ref{eq.rree}) and write
\[
\begin{split}
-\frac{1}{2s-1}\int_{\partial X}J_{1}(z,y)\>dh(y)  &  =-e^{-2\varphi(z)}\bigl
[\hat R_g(s;z,z^{\prime})-\hat R_g(1-s;z,z^{\prime})\bigr]_{z^{\prime}=z}\\
&  \qquad-\bigl[R_\tau(s;z,z^{\prime})-R_\tau(1-s;z,z^{\prime})\bigr]
_{z^{\prime}=z}
\end{split}
\]
where we note that $R_{\ast}(s)-R_{\ast}(1-s)$ has a continuous kernel and so
can be evaluated on the diagonal. In computing the finite part for $I_{1}(s)$,
the factor $e^{-2\varphi(z)}$ can be replaced by 1, since $1-e^{-2\varphi
(z)}= \mathcal{O}(\rho^{2})$. So at this stage we have
\[
\begin{split}
I_{1}(s)  &  =-{\FP_{\varepsilon\downarrow0}}\int_{\rho=\varepsilon}\bigl
[\hat R_g(s;z,z^{\prime})-\hat R_g(1-s;z,z^{\prime})\\
&  \hskip1in-R_\tau(s;z,z^{\prime})+R_\tau(1-s;z,z^{\prime})\bigr]_{z^{\prime
}=z}d\sigma_{\varepsilon}(z).
\end{split}
\]
Analysis as in Theorem 3.1 of \cite{GZ:1997} shows the integrand can by split
into two components lying in $\rho^{2}C^{0}(\bar{X})$ and $\rho^{2s}C^{\infty
}(\bar{X})$. The first contributes zero as $\varepsilon\rightarrow0$, while the
second contributes zero for $\re(s)>1$. Thus $I_{1}(s)$ vanishes for all $s$
by meromorphic continuation.

The analysis of $I_{2}(s)$ is based entirely on the fact that $(\partial_{\nu
}-s) E_{*}(s;z,y)$ vanishes to higher order as $z\to{\partial X}$ than
$E_{*}(s;z,y)$. So the conclusion that $I_{2}(s)=0$ follows exactly as in
Lemma 6.5 of \cite{PP:2000}.

The final step is to evaluate
\[
\begin{split}
I_{3}(s)  &  =-\int_{\partial X}{\FP_{\varepsilon\downarrow0}}\int_{\rho=\varepsilon
}\Bigl[e^{-2\varphi(z)}\hat E_{g}(1-s;z,y)\partial_{s}\hat E_{g}(s;z,y)\\
&  \hskip1in-E_{\tau}(1-s;z,y)\partial_{s}E_{\tau}(s;z,y)\Bigr]d\sigma_{\varepsilon}(z)dh(y).
\end{split}
\]
We will need the relation between Poisson and scattering kernels:
\[
S_{\ast}(s;y,y^{\prime})=[\rho(z)^{-s}E_{\ast}(s;z,y^{\prime})]_{z=y},
\]
for $y\neq y^{\prime}$. For convenience let $\psi_{\ast}(s;z,y)=\rho
(z)^{-s}E_{\ast}(s;z,y^{\prime})$. In the formula for $I_{3}$ above, the
$e^{-2\varphi(z)}$ may be replaced by $1$ as in the analysis of $I_{1}$.
Further more, after substituting $E_{\ast}=\rho^{s}\psi_{\ast}$, the terms
involving $\partial_{s}\rho^{s}$ are logarithmic and thus do not contribute to
the finite part. This substitution thus gives
\begin{equation}
\begin{split}
I_{3}(s)  &  =-\int_{\partial X}{\FP_{\varepsilon\downarrow0}}\int_{\rho=\varepsilon
}\Bigl[\psi(1-s;z,y)\partial_{s}\psi(s;z,y)\\
&  \hskip1in-\psi_{0}(1-s;z,y)\partial_{s}\psi_{0}(s;z,y)\Bigr]\rho \>
d\sigma_{\varepsilon}(z)dh(y).
\end{split}
\label{eq.I3calc}
\end{equation}
As pointed in \S\ref{sec.hadamard.divisor}, $S_g(s)-S_{\tau}(s)=T(s)$ is of order $2\re(s)-3$. We can
write
\begin{equation}
S_g(1-s)S^{\prime}(s)-S_{\tau}(1-s)S_{\tau}^{\prime}(s)=T(1-s)
S_g^{\prime}(s)+S_{\tau}(1-s)T^{\prime}(s), \label{eq.ssss}
\end{equation}
so this combination is of order $-2$ and has a continuous kernel. To take the
limit $\varepsilon\rightarrow0$ in $I_{3}$, we note that $\rho\>d\sigma
_{0}\rightarrow dh$ as $\rho\rightarrow0$, and the integrand in brackets in
(\ref{eq.I3calc}) approaches the kernel (\ref{eq.ssss}). Thus the limit is
\[
\begin{split}
F(s)=I_{3}(s)  &  =-\tr[S_g(1-s) S_g'(s) - S_\tau(1-s) S_\tau'(s)]\\
&  =-\frac{d}{ds}\log\det S_{g,\tau}(s).
\end{split}
\]

\section{Isopolar Surfaces, by Robert Brooks}
\label{app.brooks}

\newcommand{\overbar}{\overline}
In this appendix, we briefly review the construction of isopolar
manifolds via the Sunada method (see \cite{BP} and \cite{BGP}), and
then present a survey of various kinds of isopolar surfaces that can
be constructed in this way.  In fact, the examples are \textit{isophasal},
meaning they have identical scattering phase.  Isophasal implies
isopolar, but the eigenvalues and resonances determine the scattering
phase only up to finitely many parameters.

Let $G$ be a finite group, and let $H_1$ and $H_2$ be subgroups of $G$
satisfying the following condition, known as the Sunada condition:
\begin{equation}\label{sunada.cond}
\text{for all } g \in G, \#([g]\cap H_1) = \#([g]\cap H_2),
\end{equation}
where ``$[g]$'' denotes the conjugacy class of $g$ in $G$.

Let $L^2(G)$ denote the vector space of functions on $G$ with the left
action of $G$ given by 
$$(g^*(f))(x)= f(g^{-1}x)$$
and right action given by
$$f^g(x)=f(xg).$$

It is shown, for instance in \cite{SM} or \cite{BGP}, that (\ref{sunada.cond})
is equivalent to
\begin{equation}\label{star.cond}
(L^2(G))^{H_1}\text{ is $G$-equivalent to }
(L^2(G))^{H_2},
\end{equation}
where $(L^2(G))^{H_i}$ denotes the subspace of $L^2(G)$ which is
invariant under the left action of $H_i$.

Condition (\ref{star.cond}) is by definition equivalent to the existence of a
function 
$$c: G \to \bbR$$
such that the $G$-equivariant map
$$T: L^2(G) \to L^2(G)$$ given by
$$T(f)(x) = \sum_{g \in G} c(g)f(gx)$$ satisfies
\begin{equation}\label{tri.cond}
T \text{ is an isomorphism from }(L^2(G))^{H_1}
\text{ to } (L^2(G))^{H_2}.
\end{equation}

Note that the condition that the image of $T$ lies in $(L^2(G))^{H_2}$
is precisely that the function $c$ satisfies
$$c(gh)= c(g)\text{ for }  h \in H_2.$$

We then have:

\begin{theorem} Let $S_0$ be a surface with funnels or cusps, and let
$(G,H_1,H_2)$ satisfy (\ref{sunada.cond}) (or, equivalently, \ref{tri.cond}). 

Let $$\phi: \pi_1(S_0) \to G$$
be a homomorphism, and let $S_1$ and $S_2$ be the covers of $S_0$
whose fundamental groups are $\phi^{-1}(H_1)$ and $\phi^{-1}(H_2)$
respectively. 

Then $S_1$ and $S_2$ are isophasal.
\end{theorem}

\noindent {\bf{Sketch of proof:}} (see \cite{BGP} and \cite{BP} for
details) Let $S_{id}$ be the covering of $S_0$ whose fundamental group
is $\phi^{-1}(id)$.  We may identify $C^{\infty}(S_i)$ with
$[C^{\infty}(S_{id})]^{H_i}$. Then we have an isomorphism
$$T_S: [C^{\infty}(S_{id})]^{H_1} \to [C^{\infty}(S_{id})]^{H_2}$$ 
given by
\begin{equation}\label{tsf}
T_S(f)(x) = \sum_{g \in G} c(g)f(gx)
\end{equation}
which clearly commutes with the Laplacian. If we lift a defining
function on $S_0$ to $S_{id}$, and denote by $T_{\partial}$ the
functional on the boundary of $S_{id}$ given by the formula (\ref{tsf}), 
then $T_{\delta}$ intertwines the scattering operator on $S_1$ with
that of $S_2$. 

One way to apply this theorem is by taking $\overbar{S_0}$ to be a
closed surface (or possibly an orbifold surface) to which one can
apply the Sunada construction, and then to obtain $S_0$ from
$\overbar{S_0}$ by replacing neighborhoods of points of
$\overbar{S_0}$ with cusps or funnels, see \cite{BP} for a
discussion. In order for the singularities of $\overbar{S_0}$ to
smooth out in the coverings, we must have that the map $\phi$
satisfies a freeness condition about the singular points of
$\overbar{S_0}$. 

If we introduce a funnel (or cusp) at a regular point of
$\overbar{S_0}$, this will introduce $[G:H_i]$ funnels (or cusps) on $S_i$.
But if we introduce a funnel at an orbifold point $p$ of order $k$,
the number of funnels on $S_i$ will be $(1/k)[G:H_i]$. Note that the
freeness condition guarantees that $k$ divides $[G:H_i]$.

Using this argument, we can show:

\begin{theorem}\label{thm.iso}
\hskip4in
\begin{enumerate}
\item  There exist Riemann surfaces $S_1$ and $S_2$ which are of
genus $4$ with one funnel which are isophasal (compare \cite{BP}).

\item  There exist surfaces $S_1$ and $S_2$ of genus $3$ with one
funnel, which are conformally equivalent and carry metrics which have
constant curvature $-1$ outside of compact sets, which are isophasal
(compare \cite{BT}, \cite{BP}).

\item  There exist families of size $(c_1)k^{c_2\log{k}}$ of
mutually isophasal Riemann surfaces of genus $c_3 k$ with $c_4 k$
funnels (compare \cite{BGG}).
\end{enumerate}
\end{theorem}

One may generalize this construction in the following way: if
$\overbar{S_0}$ is an orbifold surface, we may take $\phi$ to be any
homomorphism 
$$\phi: \pi_1(\overbar{S_0}) \to G,$$
which need not satisfy the freeness condition on all the orbifold
points of $\overbar{S_0}$. We may then introduce funnels at points of
$\overbar{S_0}$ as long as we introduce funnels at all the orbifold
points which do not satisfy the freeness condition. 

With this idea, we can show:

\begin{theorem}
\hskip4in
\begin{enumerate}
\item  There exist two surfaces of genus $2$ with four funnels
which are isophasal.

\item  There exist two Riemann surfaces of genus $3$ with three
funnels which are isophasal.

\item  There exist two surfaces of genus $0$ with eight funnels
which are isophasal.

\item There exists two Riemann surfaces of genus $0$ with sixteen
funnels which are isophasal.

\end{enumerate}
\end{theorem}

Details of these examples will appear elsewhere.

\end{document}